\documentclass{article}
\usepackage{amsmath,amsfonts,amsthm,epsfig,amssymb}

\newtheorem{theorem}{Theorem}
\newtheorem{lemma}[theorem]{Lemma}

\theoremstyle{definition}
\newtheorem{remark}[theorem]{Remark}

\renewcommand\Pr{{\mathop{\mathbb P{}}\nolimits}}
\renewcommand\phi{\varphi}

\newcommand\op{{o_{\mathrm{p}}}}
\newcommand\E{\operatorname{\mathbb E{}}}

\newcommand\Bi{\operatorname{Bi}}
\newcommand\ceil[1]{{\lceil #1 \rceil}}

\newcommand\eps{\varepsilon}

\newcommand\cc{{\mathrm c}}

\newcommand\Gkl{G^{k,\ell}}
\newcommand\Gokl{{\widetilde G}^{k,\ell}}
\newcommand\U{{\mathcal U}}
\newcommand\D{{\mathcal D}}
\newcommand\bp{{\mathfrak X}}
\newcommand\bpm{{\mathfrak X}'}
\newcommand\la{\lambda}
\newcommand\noproof{\hfill$\Box$}
\newcommand\nv{{\nu}}
\newcommand\bb[1]{\bigl(#1\bigr)}
\newcommand\ore[1]{\protect\overrightarrow{#1}}% protect needed in figure captions
\newcommand\Gd{\ore{G}}
\newcommand\Gdkl{\ore{G}^{k,\ell}}
\newcommand\Hd{\ore{H}}

\newcommand\Kdk{\ore{K_k}}
\newcommand\Kdl{\ore{K_\ell}}
\newcommand\mud{\ore{\mu}}
\newcommand\aut{{\mathrm{aut}}}
\newcommand\F{{\mathcal F}}
\newcommand\tF{{\widetilde{\mathcal F}}}

\newcommand\evA{{\mathcal A}}
\newcommand\evB{{\mathcal B}}
\newcommand\evC{{\mathcal C}}
\newcommand\evD{{\mathcal D}}
\newcommand\evE{{\mathcal E}}
\newcommand\evF{{\mathcal F}}
\newcommand\evU{{\mathcal U}}
\newcommand\evIt{{\mathcal I}_t}
\newcommand\evQ{{\mathcal Q}}
\newcommand\tD{{\widetilde \evD}}
\newcommand\tU{{\widetilde \evU}}
\newcommand\eva{{\mathcal U}}% the increasing event in the FK argument
\newcommand\eve{{\mathcal E}}% from same place

\begin{document}
\title{Clique percolation}
\date{September 8, 2008}

\author{B\'ela Bollob\'as\thanks{Department of Pure Mathematics
and Mathematical Statistics, University of Cambridge, Cambridge CB3 0WB, UK}
\thanks{University of Memphis, Memphis TN 38152, USA}
\thanks{Research supported in part by NSF grants DMS-0505550,
 CNS-0721983 and CCF-0728928, and ARO grant W911NF-06-1-0076}
\and Oliver Riordan
\thanks{Mathematical Institute, University of Oxford, 24--29 St Giles', Oxford OX1 3LB, UK}}
\maketitle

\begin{abstract}
Der{\'e}nyi, Palla and Vicsek introduced the following dependent
percolation model, in the context of finding communities in networks.
Starting with a random graph $G$ generated by some rule, form an auxiliary
graph $G'$ whose vertices are the $k$-cliques of $G$, in which two 
vertices are  joined if the corresponding cliques share $k-1$ vertices.
They considered in particular the case where $G=G(n,p)$, and found
heuristically the threshold function $p=p(n)$ above which
a giant component appears in $G'$.
Here we give a rigorous proof of this result, as well as many extensions.
The model turns out to be very interesting due to the essential
global dependence present in $G'$.
\end{abstract}

\section{Cliques sharing vertices}\label{sec_cv}

Fix $k\ge 2$ and $1\le \ell\le k-1$. Given a graph
$G$, let $\Gkl$ be the graph whose vertex set is the set of all copies of $K_k$
in $G$, in which two vertices are adjacent if the corresponding copies
of $K_k$ share at least $\ell$ vertices.
Starting from a random graph $G=G(n,p)$, our aim is to study percolation
in the corresponding graph $\Gkl_p$, i.e., to find for which values
of $p$ there is a `giant' component
in $\Gkl_p$, containing a positive fraction of the vertices of $\Gkl_p$.

For $\ell=k-1$, this question was
proposed by Der{\'e}nyi, Palla and Vicsek~\cite{DPV},
motivated by the study of `communities' in real-world networks,
but independent of the motivation, we consider it to be an extremely
natural question in the theory of random graphs. Indeed, it is perhaps
the most natural example of dependent percolation arising out
of the model $G(n,p)$.

As we shall see in a moment, it is not too hard to guess the answer;
simple heuristic derivations based on the local analysis
of $\Gkl_p$ were given in~\cite{DPV} and by
Palla, Der{\'e}nyi and Vicsek~\cite{PDV}.
(For a survey of related work see~\cite{PAFPDV}.)
Note, however, that $\Gkl_p$ may well have many more than $n^2$ edges, so $\Gkl_p$
is not well approximated by a graph with independence between different edges: there
is simply not enough information in $G(n,p)$.
Thus it is not surprising that it requires
significant work to pass from local information about $\Gkl_p$
to global information about the giant component.
Nonetheless, it turns out
to be possible to find exactly the threshold for percolation,
for all fixed $k$ and $\ell$.

Given $0<p=p(n)<1$, let
\begin{equation}\label{muform}
 \mu = \mu(n,p) = \left(\binom{k}{\ell}-1\right)\binom{n}{k-\ell} p^{\binom{k}{2}-\binom{\ell}{2}},
\end{equation}
so $\mu$ is $\binom{k}{\ell}-1$ times the expected number of $K_k$s
containing a given copy of $K_{\ell}$. Intuitively,
this corresponds to the average number of new $K_{\ell}$s reached
in one step from a given $K_{\ell}$, so we expect percolation
if and only if $\mu>1$.
Since $\binom{k}{2}-\binom{\ell}{2} = \ell(k-\ell)+\binom{k-\ell}{2}
= (k-\ell)(k+\ell-1)/2$, we have
$\mu=\Theta(1)$ if and only if
\begin{equation}\label{theta}
 p = \Theta\left(n^{-\frac{2}{k+\ell-1}}\right);
\end{equation}
we shall focus our attention on $p$ in this range.

In addition to finding the threshold for percolation, we shall
also describe the asymptotic proportion
of $K_k$s in the giant component in terms of the survival probability
of a certain branching process.
Set $M=\binom{k}{\ell}-1$.
Given $\la>0$, let $Z_\la$ have a Poisson distribution with
mean $\la/M$. Let $\bp(\la)=(X_t)_{t=0}^\infty$
be the Galton--Watson branching process which starts with a single
particle in $X_0$, in which each particle in $X_t$ has children
in $X_{t+1}$ independently of the other particles and of
the history, and in which the distribution of the number of children
of a given particle is given by $MZ_\la$.
Let $\rho=\rho(\la)$ denote the probability that $\bp(\la)$ does not die out.
Then a simple calculation shows that $\rho$
satisfies the equation
\[
 \rho= 1-\exp\left(-(\la/M)(1-(1-\rho)^M)\right).
\]
From standard branching process results,
$\rho$ is the largest solution to this equation,
$\rho(\la)$ is a continuous function of $\la$, and $\rho(\la)>0$
if and only if $\la$, the expected number of children of each
particle, is strictly greater than $1$.

Let $\bpm(\la)$ denote the union of $\binom{k}{\ell}$ independent
copies of the branching process $\bp(\la)$ described above, and
let $\sigma=\sigma(\la)$ denote the survival probability of $\bpm(\la)$,
so $\sigma=1-(1-\rho)^{\binom{k}{\ell}}$.
Our main result is that when $\mu=\Theta(1)$, the largest component
of $\Gkl_p$ contains whp a fraction $\sigma(\mu)+o(1)$
of the vertices of $\Gkl_p$, where $\mu$ is defined by \eqref{muform}.
Here, as usual, an event holds {\em with high probability}, or whp, if its
probability tends to $1$ as $n\to\infty$.

Let $\nv=\binom{n}{k}p^{\binom{k}{2}}$ denote
the expected number of copies of $K_k$ in $G(n,p)$,
i.e., the expected number of vertices of $\Gkl_p$.
Let us write $C_i(G)$ for the number of vertices in
the $i$th largest component of a graph $G$.

\begin{theorem}\label{th1}
Fix $1\le \ell<k$,
and let $p=p(n)$ be chosen so that $\mu=\Theta(1)$, where $\mu$ is defined by \eqref{muform}.
Then, for any $\eps>0$, whp we have
\[
 (\sigma(\mu)-\eps) \nv \le C_1(\Gkl_p) \le
 (\sigma(\mu)+\eps) \nv
\]
and $C_2(\Gkl_p)\le \eps \nv$.
\end{theorem}
It is well known that $|\Gkl_p|$ is concentrated
around its mean $\nv$ whenever $\nv\to\infty$,
so Theorem~\ref{th1} simply says
that the largest component of $\Gkl_p$ contains a fraction $\sigma(\mu)+o(1)$
of the vertices whp. The extension to the case where $\mu\to 0$
or $\mu\to\infty$ is essentially trivial, and will be discussed
in Subsection~\ref{ss_ext}.

We shall prove Theorem~\ref{th1} in two stages, considering
the subcritical case in
the next subsection, and the supercritical case in Subsection~\ref{ss_lb}.
Very roughly speaking, to handle the subcritical case (and to prove
the upper bound on the giant component in the supercritical case) we shall
show {\em approximate} domination of a suitable component exploration
in $\Gkl_p$ by the branching process $\bp'(\la)$, $\la=(1+\eps)\mu$.
Due to the dependence in the model, we have to be very careful exactly
how we explore $\Gkl_p$ to make this argument work. For the upper bound
we first show (by approximate local coupling with the branching process)
that roughly the right number of vertices are in large components, even
if $p$ is reduced slightly, i.e., even if we omit some edges.
Then we use a multi-round `sprinkling' argument, putting back the omitted
edges in several rounds, and showing that it is very likely that the sprinkled
edges join these large components. The details of both arguments turn
out to be less simple that one might like.

\subsection{The subcritical case}\label{ss_ub}

We shall start by considering the subcritical case, proving
the following much stronger form of Theorem~\ref{th1} in
this case.

\begin{theorem}\label{th1l}
Let $1\le \ell\le k-1$ and $\eps>0$ be given.
There is a constant $C=C(k,\ell,\eps)$ such that,
if $p=p(n)$ is chosen so that $\mu\le 1-\eps$
for all large enough $n$, then $C_1(\Gkl_p)\le C\log n$ whp.
\end{theorem}

\begin{proof}
Since the event $C_1(\Gkl)> C\log n$, considered as a property
of the underlying graph $G$, is an increasing event, we may assume
without loss of generality that $\mu=1-\eps$
for every $n$. Thus \eqref{theta} holds.

Fixing a set $V_0$ of $k$ vertices of $G=G(n,p)$, we shall show that,
given that $V_0$ forms a complete graph in $G$,
the probability that the corresponding component $C(V_0)$ of $\Gkl_p$
has size more than $C\log n$ is at most $n^{-k-1}$, provided
$C$ is large enough. Since the probability that $V_0$
forms a complete graph in $G$ is $p^{\binom{k}{2}}$, while
there are $\binom{n}{k}$ possibilities for $V_0$,
it then follows that
$\Pr(C_1(\Gkl_p)\ge C\log n)\le \binom{n}{k}p^{\binom{k}{2}}n^{-k-1}=o(1)$.

From now on we condition on $V_0$ forming a $K_k$ in $G=G(n,p)$.
The strategy is to show domination of a natural component exploration
process by the branching process described earlier.
We shall show essentially that the average number of new $K_\ell$s
reached from a given $K_\ell$ in $G$ via $K_k$s in $G$
is at most $\mu+o(1)$,
though there will be some complications.

In outline, our exploration of the component $C(V_0)\subset \Gkl_p$
proceeds as follows.
At each stage we have a set $V_t$ of {\em reached} vertices of $G$, starting
with $V_0$; we
also keep track of a set $E$ of {\em reached} edges, initially
the edges spanned by $V_0$. 
At the end of stage $t$ of our exploration, $E$ will consist
of all edges of $G[V_t]$.
Within $V_t$, every $K_{\ell}$ is labelled as either
`tested' or `untested'.
We start with all $\binom{k}{\ell}$
$K_{\ell}$s in $V_0$ marked as untested.
The exploration stops when there are no untested $K_{\ell}$s.

As long as there are untested $K_{\ell}$s, we proceed as follows.
Pick one, $S$, say.
One by one, test each set $K$ of $k$ vertices with $S\subset K\not\subset V_t$
to see whether all edges induced by $K$ are present in $G$.
If so, we add any new vertices to $V_t$, i.e., we set $V_{t+1}=V_t\cup V(K)$.
We now add all edges of $K$ not present in $V_t$ to $E$;
we call these edges {\em regular}.
Any new $K_{\ell}s$ formed in $E$ are marked as untested.
Note that any such $K_{\ell}$ must contain at least one vertex
of $V_{t+1}\setminus V_t$, and hence must lie entirely inside $K$.

Next, we test all edges
between $V_t$ and $V(K)\setminus V_t$ to see if they are present
in $G$, adding any edge found to $E$,
and marking any new $K_{\ell}$s
formed as untested. Edges added during this operation
are called  {\em exceptional}.
At this point, we have revealed the entire subgraph of $G$ induced
by $V_{t+1}$, i.e., we have $E=E(G[V_{t+1}])$.
We then continue through our list of possible sets $K$ containing $S$, omitting
any set $K$ contained in the now larger set $V_{t+1}$.
Once we have considered all possible $K\supset S$, we mark
$S$ as tested, and continue to another untested $K_{\ell}$, 
if there is one.

The algorithm described above can be broken down into a sequence
of steps of the following form.
At the $i$th step, we test whether all edges
in a certain set $A_i$ are present in $G=G(n,p)$; the future
path of the exploration depends only on the answer to this question,
not on which particular edges are missing if the answer is no.
Although this is wasteful from an algorithmic point of view, it
is essential for the analysis.
We write $\evA_i$ for the event $A_i\subset E(G)$.
After $i$ steps, we will have `uncovered' a set $E_i$
of edges (called $E$ above). The set $E_i$ consists of the
edges spanned by $V_0$ together with the union of those sets $A_j$ for
which $\evA_j$ holds.

The event that the algorithm reaches a particular state, i.e.,
receives a certain sequence of answers to the first $i$ questions,
is of the form $\U\cap \D$, where $\U=\{E_i\subset E(G)\}$
is an up-set, and $\D$ is a down-set,
formed by the intersections of various $\evA_j^\cc$.
The key point is that $\U$ is a {\em principal} up-set,
so $\U\cap\D$ may be regarded as a down-set $\D'$ in the product
probability space $\Omega'=\{0,1\}^{E(K_n)\setminus E_i}$
with the appropriate measure.
Hence, for any $A_{i+1}$ disjoint from $E_i$, the conditional
probability that $\evA_{i+1}$ holds given the current state of the algorithm
is
\begin{equation}\label{dd}
 \Pr(\evA_{i+1}\mid \U\cap \D) = \Pr(\evA_{i+1}'\mid \D') \le \Pr(\evA_{i+1}') = p^{|A_{i+1}|},
\end{equation}
where $\evA_{i+1}'$ is the event in $\Omega'$ corresponding to $\evA_{i+1}$,
and the inequality follows from Harris's Lemma~\cite{Harris} applied in $\Omega'$.

Let us write $X_i$ for the number of new $K_{\ell}$s found as a result
of adding regular edges when testing the $i$th $K_{\ell}$, $S_i$, say;
we shall deal with exceptional edges separately in a moment.
Recall that we add regular edges when we find a new $K_k$ with at most
$k-\ell$ and at least $1$ vertex outside the current vertex set $V_t$.

Let $\eta>0$ be a constant such that $(1+\eta)\mu\le (1-\eps/2)$.
When testing $S_i$, there are at most $\binom{n}{k-\ell}$
possibilities for new $K_k$s with $k-\ell$ vertices outside
the current $V_t$. Given the history, by \eqref{dd} each such $K_k$
is present with probability at most
$p^{\binom{k}{2}-\binom{\ell}{2}}$,
so the number of such $K_k$s we find is stochastically dominated
by the Binomial distribution
$\Bi\left(\binom{n}{k-\ell},p^{\binom{k}{2}-\binom{\ell}{2}}\right)$,
and hence, for $n$ large, by a Poisson distribution with mean
$(1+\eta/2)\binom{n}{k-\ell}p^{\binom{k}{2}-\binom{\ell}{2}}$.
[ Here we use the fact that a Poisson distribution with
mean $-N\log(1-\pi)$ dominates a Binomial $\Bi(N,\pi)$,
which, as pointed out to us by Svante Janson,
follows immediately from the same statement for $N=1$. ]

For $1\le j\le k-\ell-1$, we may also find new $K_k$s
containing $S_i$
together with $j$ other vertices of the current set $V_t$,
and hence with only $k-\ell-j$ vertices outside $V_t$.
Assuming $|V_t|\le k(\log n)^{100k^3} \le (\log n)^{101k^3}$, say, 
the number of possibilities for a fixed $j$ is crudely
at most $(\log n)^{101k^3j} n^{k-\ell-j}$, and each of these
tests succeeds with probability at most $p^{\binom{k}{2}-\binom{\ell+j}{2}}$.
A simple calculation shows that
$n^{k-\ell-j}p^{\binom{k}{2}-\binom{\ell+j}{2}}$ is at most
$n^{-\delta}$ for some $\delta>0$, so the expected number
of $K_k$s of this type is at most $n^{-\delta/2}$, say.
Moreover, the distribution of the number found is stochastically dominated
by a Poisson distribution with mean $2n^{-\delta/2}$.

Each $K_k$ we find consisting of $k-\ell-j$ new vertices and $\ell+j$
old vertices, $j\ge 0$, generates $\binom{k}{\ell}-\binom{\ell+j}{\ell}
\le M$ new $K_{\ell}$s, where $M=\binom{k}{\ell}-1$.
It follows that, given the history, the conditional distribution
of $X_i/M$ is stochastically dominated by a Poisson distribution with mean
\begin{equation}\label{bpdom}
 (1+\eta/2)\binom{n}{k-\ell}p^{\binom{k}{2}-\binom{\ell}{2}} + 2n^{-\delta/2}
 = (1+\eta/2)\mu/M +o(1),
\end{equation}
which is at most $(1+\eta)\mu/M<(1-\eps/2)/M$ if $n$ is large enough.

Turning to exceptional edges, we claim that the $j$th exceptional edge
added creates at most $\binom{k-1+j}{\ell-1}$ new $K_{\ell}$s;
all we shall use about this bound is that it
depends only on $j$, $k$ and $\ell$, not on $n$.
Indeed, we add exceptional edges immediately after adding
a $K_k$ that includes a certain set $N$ of new vertices.
At this point, the degree in $E$ (the uncovered edges)
of every vertex in $N$ is exactly $k-1$. We now add one or
more exceptional edges joining $N$ to $V_t$.
Any such edge $e$ has one end, $x$, say, in $N$.
If $e$ is the $j$th exceptional edge in total,
then just after adding $e$ the vertex $x$ has degree at most $k-1+j$.
Any new $K_{\ell}$s involving $e$ consist of $x$
together with $\ell-1$ neighbours of $x$,
so there are at most $\binom{k-1+j}{\ell-1}$
such $K_{\ell}$s.

Assuming $|V_t|\le k(\log n)^{100k^3}$, the number of potential exceptional
edges associated to a new $K_k$ is at most $(k-\ell)|V_t|=O^\star(1)$,
where, as usual, $g_1(n)=O^\star(g_2(n))$ means that
there is a constant $a$ such that $g_1(n)=O(g_2(n)(\log n)^a)$.
It follows that, for fixed $r$, the probability that we find at least $r$
such edges at a given step is $O^\star(p^r)$.
Furthermore,
the probability that we find $j$ exceptional edges
in total during the first $(\log n)^{100k^3}$ steps is $O^\star(p^j)$,
since there
are $O^\star(1)$ possibilities for the set of at most $j$ steps
at which we might find them.
Let us choose a constant $J$ so that $p^J\le n^{-100k^3}$
(here, $p^J\le n^{-k-2}$ would do; the stronger bound
is useful later),
and let $\evB$ be the `bad' event that we find more than $J$ exceptional
edges in the first $(\log n)^{100k^3}$ steps. Then
we have $\Pr(\evB)=O^\star(p^J)=O^\star(n^{-100k^3})=o(n^{-99k^3})$.

As long as $\evB$ does not hold, we create
at most $J'=\sum_{j\le J}\binom{k-1+j}{\ell-1}=O(1)$
new $K_{\ell}$s when adding exceptional edges
in the first $(\log n)^{100k^3}$ steps;
let us note for later that we also create at most $\sum_{j\le J}\binom{k-1+j}{k-1}$
$K_k$s when adding exceptional edges.
We view our exploration as a set of
branching processes: we start one process for each of the initial
$K_{\ell}$s. Whenever we add a $K_{\ell}$ in the normal
way, we view it as a child of the $K_{\ell}$ we were testing.
When we add a $K_{\ell}$ as a result of adding an exceptional
edge, we view it as the root of a new process.
As long as $|V_t|\le k(\log n)^{100k^3}$ holds,
from~\eqref{bpdom} the branching processes
we construct are stochastically dominated
by independent copies $\bp_i$ of the Galton--Watson process $\bp(\la)$
described earlier, where $\la=(1+\eta)\mu < (1-\eps/2)$.
If $\evB$ does not hold, then we start in total at most $J''=\binom{k}{\ell}+J'=O(1)$
processes in the first $(\log n)^{100k^3}$ steps.

Recall that the offspring distribution in $\bp(\la)$ is given by $M Z_\la$,
where $Z_\la$ has a Poisson distribution with mean $\la/M$, so $\E(MZ_\la)=\la$.
Here, $\la=(1+\eta)\mu < 1-\eps/2$.
Since $M Z_\la$ has an exponential upper tail, it follows 
from standard branching process results that there is a constant $a>0$
such that the probability that the total size
of $\bp_i$ exceeds $m+1$ is at most $\exp(-am)$ for any $m\ge 0$.
Taking $C$ large enough, it follows that with probability
$1-o(n^{-k-1})$, each of $\bp_1,\ldots,\bp_{J''}$
has size at most $(C/J'')\log n$. If this event holds and $\evB$
does not hold, then our exploration dies having reached a total
of at most $C\log n$ vertices.
Hence, the probability that $C(V_0)$ contains more than $C\log n\le (\log n)^{100k^3}$
vertices of $G=G(n,p)$ is $o(n^{-k-1})+o(n^{-99k^3})=o(n^{-k-1})$.

At this point the proof of Theorem~\ref{th1l} is almost complete:
we have shown that whp, any component
of $\Gkl_p$ involves $K_k$s meeting at most $C\log n$
vertices of $G=G(n,p)$.
To complete the proof,
it is an easy exercise to show that if $p\le n^{-\delta}$ for
some $\delta>0$, then whp any $C\log n$ vertices of $G(n,p)$
span at most $C'\log n$ copies of $K_k$, for some constant $C'$.
Alternatively, note that
the number of $K_k$s found involving new vertices is at most the
final number of vertices reached, while all other $K_k$s are
formed by the addition of exceptional edges, and if $\evB$
does not hold, then, arguing as for the bound on the number
of $K_\ell$s formed by adding exceptional edges,
the number of $K_k$s so formed is bounded by a constant.
\end{proof}

In the proof above, subcriticality only came in at the end,
where we used it to show that the branching processes $\bp_i$
were very likely to die; in the supercritical case,
the proof gives a domination result that we shall state
in a moment. For this, the order in which we test the
$K_{\ell}$s matters -- we proceed in rounds, in round
$0$ testing the $\binom{k}{\ell}$ initial $K_{\ell}$s,
and then in round $i\ge 1$ testing all $K_{\ell}$s
created during round $i$.
Let $H=H(\Gkl_p)$ be the bipartite incidence graph corresponding
to $\Gkl_p$: the vertex classes are $V_1$, the set
of all $K_k$s in $\Gkl_p$, and $V_2$, the set of all
$K_{\ell}$s. Two vertices are joined if one of the
corresponding complete graphs is contained in the other.
Given a vertex $v_0\in V_1$ of $H$, let $N_i=N_i(v_0)$
denote the number of $K_{\ell}$s whose graph distance in $H$
from $v_0$ is at most $2i+1$.
If $v_0$ is the vertex of $H$ corresponding to 
the complete subgraph on $V_0$, then after $i$ rounds
of the above algorithm we have certainly reached
all $N_i$ $K_{\ell}$s within distance $2i+1$
of $v_0$.

The domination argument in the proof of Theorem~\ref{th1l}
thus also proves the lemma below, in which $J''$ is a constant
depending only on $k$ and $\ell$, $\bp_1,\ldots,\bp_{J''}$
are independent copies of our Galton--Watson branching process
as above, and $M_{\le t}(\bp_1,\ldots,\bp_{J''})$
denotes the total number of particles
in the first $t$ generations of $\bp_1,\ldots,\bp_{J''}$.

\begin{lemma}\label{ub}
Let $\eta>0$ be fixed,
let $p=p(n)$ satisfy \eqref{theta}, and let $V_0$ be a fixed set
of $k$ vertices of $G=G(n,p)$. Condition on $V_0$ spanning a complete
graph in $G$, and let $v_0$ be the corresponding vertex of $H$.
Then we may couple the random sequence $N_1,N_2,\ldots$
with $J''$ independent copies $\bp_i$ of $\bp((1+\eta)\mu)$
so that with probability $1-o(n^{-99k^3})$ we have
$N_t\le M_{\le t}=M_{\le t}(\bp_1,\ldots,\bp_{J''})$
for all $t$ such that $M_{\le t}\le (\log n)^{100k^3}$.
\noproof
\end{lemma}

We finish this subsection by presenting a consequence 
of a much simpler version of the domination argument above.
If we are prepared to accept a larger
error probability, we may abandon the coupling
the first time an exceptional edge appears. As shown
above, the probability that we find any exceptional edges
within $O^\star(1)$ steps is at most $n^{-\delta}$ for some $\delta>0$.
Abandoning our coupling if this happens, we need only
consider the original $\binom{k}{\ell}$ branching processes, one for each copy
of $K_{\ell}$ in $V_0$.
In other words, we may compare our neighbourhood exploration
process with the branching process $\bpm(\la)$, $\la=(1+\eta)\mu$, which
starts with $\binom{k}{\ell}$ particles in generation $0$,
and in which, as in $\bp(\la)$, the offspring distribution
for each particle is given by $M$ times a Poisson distribution
with mean $\la/M$.

\begin{lemma}\label{ubc}
Let $\eta>0$ be fixed,
let $p=p(n)$ satisfy \eqref{theta}, and let $V_0$ be a fixed set
of $k$ vertices of $G=G(n,p)$.
Condition on $V_0$ spanning a complete
graph in $G$, and let $v_0$ be the corresponding vertex of $H$.
Then there is a constant $\delta>0$ such that
we may couple the random sequence $N_1,N_2,\ldots$
with $\bpm=\bpm((1+\eta)\mu)$
so that, with probability at least $1-n^{-\delta}$, we have
$N_t\le M_t$
for all $t$ such that $M_{\le t}\le (\log n)^{100k^3}$,
where $M_t$ is the number of particles in generation $t$ of $\bpm$,
and $M_{\le t}=M_0+M_1+\cdots+M_t$.
\end{lemma}

In the next subsection we shall show that when $\mu>1$,
the graph $\Gkl_p$ does contain a giant component, and moreover that
this giant component is of about the right size; Lemma~\ref{ubc}
will essentially give us the upper bound, but we have
to work a lot more for the lower bound.

\subsection{The supercritical case}\label{ss_lb}

Recall that
$|\Gkl_p|$, the number of $K_k$s in $G(n,p)$, 
is certainly concentrated about its mean $\nv=\binom{n}{k}p^{\binom{k}{2}}$.
For the moment, we concentrate on the case where \eqref{theta}
holds; we return to larger $p$ later.

One bound in Theorem~\ref{th1} is easy,
at least in expectation: Lemma~\ref{ubc}
gives an upper bound on the expected size of the giant component.
In fact, it gives much more, namely an upper bound on the expected
number of vertices in `large' components. It is convenient
to measure the size of a component by the number of $K_{\ell}$s
rather than the number of $K_k$s.

Let $N_{\ge a}(\Gkl_p)$ denote the number of vertices of $\Gkl_p$
whose component in the bipartite graph $H$ contains at least
$a$ vertices of $V_2$, i.e., at least $a$ copies of $K_{\ell}$.

\begin{lemma}\label{ll1}
Let $p=p(n)$ be chosen so that $\mu$ is constant, and let
$\eps>0$ be fixed. For any $\omega=\omega(n)$ tending to infinity
we have $\E(N_{\ge \omega}(\Gkl_p)) \le (\sigma(\mu)+\eps)\nv$
if $n$ is large enough.
\end{lemma}
\begin{proof}
We may assume without loss of generality that $\omega\le\log n$.
From standard branching process results, for any fixed $\la$,
the probability that $\bpm(\la)$ contains at least $a$
particles but does not survive forever tends to $0$ as $a\to\infty$.
Thus, $\Pr(|\bpm(\la)|\ge \omega)=\sigma(\la)+o(1)$.

Fix a set $V_0$ of $k$ vertices of $G=G(n,p)$,
and condition on $V_0$ forming a $K_k$ in $G$, which we denote $v_0$.
Let $\la=(1+\eta)\mu$ where, for the moment, $\eta>0$ is constant.
Since $\omega\le \log n$, Lemma~\ref{ubc} tells us that the probability
$\pi$ that the component of $v_0$ in $H$ contains at least $\omega$
$K_{\ell}$s is at most
$\Pr(|\bpm((1+\eta)\mu)|\ge \omega)+o(1)=\sigma((1+\eta)\mu)+o(1)$.
Letting $\eta\to 0$ and using continuity of $\sigma$, it follows that
$\pi\le \sigma(\mu)+o(1)$ as $n\to\infty$.
Since $\E(N_{\ge\omega})$ is simply $\pi\binom{n}{k}p^{\binom{k}{2}}$,
this proves the lemma.
\end{proof}

As in Bollob\'as, Janson and Riordan~\cite{BJR}, for example,
a simple variant of Lemma~\ref{ubc} also gives us a second moment bound.

\begin{lemma}\label{ll2}
Let $p=p(n)$ be chosen so that $\mu$ is constant, and let
$\eps>0$ be fixed. For any $\omega=\omega(n)$ tending to infinity
we have $\E\bb{N_{\ge \omega}(\Gkl_p)^2} \le (\sigma(\mu)^2+\eps)
\nv^2$.
\end{lemma}
\begin{proof}
The expected number of pairs of overlapping $K_k$s in $G=G(n,p)$
is
\[
 \sum_{i=1}^{k-1}\binom{n}{k}\binom{k}{i}\binom{n-k}{k-i}p^{2\binom{k}{2}-\binom{i}{2}},
\]
which, by a standard calculation, is $o(\nv^2)$.
Hence, it suffices to bound the expected number of pairs of vertex disjoint
$K_k$s each in a `large' component. We may do so as in the proof of
Lemma~\ref{ll1}, using a variant of Lemma~\ref{ubc} in which
we start with two disjoint $K_k$s, and explore from each separately,
abandoning each exploration if it reaches size at least $\log n$,
and abandoning both if they meet, an event of probability $o(1)$.
\end{proof}

Let us turn to our proof of the heart of Theorem~\ref{th1},
namely the lower bound.
In proving this we may assume that $\mu>1$ is constant.
We start with a series of simple lemmas.

Let $V_0$ be a set of $k$ vertices of $G=G(n,p)$,
and let $\evA=\evA(V_0)$ be the event that $V_0$ spans a $K_k$ in $G$.
Let $\evQ=\evQ(V_0)$ be the event that $\Gkl_p$ contains a tree $T$
with  $\ceil{(\log n)^{5k^3}}$ vertices,
one of which, $v_0$, is the clique corresponding to $V_0$,
with the following additional property:
ordering the vertices of $T$ so that the distance from
the root $v_0$ is increasing,
each corresponding $K_k$ meets the union of all earlier $K_k$s
in exactly $\ell$ vertices. Equivalently, the union of the cliques
in $T$ contains exactly $k+(k-\ell)(|T|-1)$ vertices of $G$.

Recall that $\mu=\mu(n,p)$ is defined by \eqref{muform}.

\begin{lemma}\label{el}
Fix $\eps>0$, and let $p=\Theta(n^{-\frac{2}{k+\ell-1}})$ be chosen 
so that $\mu$ is a constant greater than $1$.
If $n$ is large enough, then $\Pr(\evQ\mid \evA)\ge \sigma(\mu)-\eps$.
\end{lemma}

\begin{proof}
Throughout we condition on $\evA=\evA(V_0)$, writing $v_0$ for the corresponding
vertex of $\Gkl_p$. We start by marking all $\binom{k}{\ell}$
copies of $K_{\ell}$ in $V_0$ as untested;
we shall then explore part of the component
of $\Gkl_p$ containing the vertex $v_0$ corresponding to $V_0$.
At the $i$th step in our exploration,
we consider an untested copy $S_i$ of $K_{\ell}$, and test for the presence
of certain $K_k$s consisting of $S_i$ plus exactly $k-\ell$ `new' vertices
not so far reached in our exploration. For each such $K_k$ we find,
we mark the $M=\binom{k}{\ell}-1$ new $K_{\ell}$s created as untested;
having found all such $K_k$s,
we mark $S_i$ as tested. We abandon our exploration
if there is no untested $S_i$ left, or if we reach more than $(\log n)^{5k^3}$
$K_k$s. Note that the total number of vertices reached
is exactly $|V_0|$ plus $k-\ell$ times the number of $K_k$s found,
so if we find more than $(\log n)^{5k^3}$ $K_k$s, then $\evQ(V_0)$ holds.

The exploration above corresponds to the construction of $\binom{k}{\ell}$
random rooted trees
whose vertices are the $S_i$, in which the children of $S_i$
are the new $K_{\ell}$s created when testing $S_i$. The number
of children of $S_i$ is $MX_i$, where $X_i$ is the number of $K_k$s we find
when testing $S_i$. Let $0<\eta<1$ be a constant to be chosen
later. Let $Z_1,Z_2,\ldots$
be a sequence of iid Poisson random variables with mean
$(1-\eta)\mu/M<\mu/M$. Our aim is to show that as long as we have
found at most $(\log n)^{5k^3}$ copies of $K_k$ in total, the conditional
distribution of $X_i$ given the history may be coupled with $Z_i$ so that
$X_i\ge Z_i$ holds with probability $1-o(n^{-\delta})$, for some
$\delta>0$. The Galton--Watson branching process $\bpm((1-\eta)\mu)$
defined by $Z_1,Z_2,\ldots$ is supercritical,
and so survives forever with probability $\sigma((1-\eta)\mu)$.
It then follows that $\evQ(V_0)$ holds
with probability at least $\sigma((1-\eta)\mu)-o(1)$.
Using continuity of $\sigma$ and choosing $\eta$ small enough, 
the conclusion of the lemma follows.

In order to establish the coupling above, we must be a little careful
with the details of our exploration. At step $i$, before
testing $S_i$, we will have a certain set $V_i$ of reached
vertices, consisting of all vertices of all $K_k$s found so far,
and a certain set $D_i\supset V_i$ of `dirty' vertices. The remaining
vertices are `clean'; we write $C_i$ for the set of these vertices.
At the start, $V_0$ is our initial set of $k$ vertices, while $D_0=V_0$
and $C_0=V(G)\setminus V_0$.

We test $S_i$ as follows: for each $v\in C_i$, let $\evE_{v,i}$ be the
event that all $\ell$ possible edges joining $v$ to $S_i$
are present in $G=G(n,p)$. First, for every vertex
$v\in C_i$, we test whether $\evE_{v,i}$ holds, writing $W_i$ for the
set of $v\in C_i$ for which $\evE_{v,i}$ does hold.
We then look for copies of $K_{k-\ell}$ inside $G[W_i]$,
writing $N_i$ for the maximum number of vertex disjoint copies.
Taking a particular set of $N_i$ disjoint copies, we then
add each of the corresponding $K_k$s to our component,
defining $V_{i+1}$ appropriately.
We then set $D_{i+1}=D_i\cup W_i$, and $C_{i+1}=V(G)\setminus D_{i+1}$.

The structure of the algorithm guarantees the following:
given the state at time $i$, all we know about the edges
between $V_i$ and $C_i$ is that certain sets of $\ell$
edges are not all present: more precisely, we know
exactly that none of the events $\evE_{v,j}$ holds, for $v\in C_i$ and $j<i$.
Let $n_i=|W_i|$, a random variable.
Having found $W_i$, it follows that the edges within $W_i$
are untested, so each is present with its unconditional probability,
and $G[W_i]$ has the distribution of the random graph $G(n_i,p)$.
Let $\eta'>0$ be a very small constant
to be chosen below. Let $\evE_i$ be the event
that $n_i\ge (1-2\eta')np^{\ell}$.

We shall show in a moment that $\evE_i$ holds with very high (conditional)
probability, given the history; first, let us see
how this enables us to complete the proof.
If $\evE_i$ does hold, then the conditional expected number
of $K_{k-\ell}$s in $G[W_i]$ is exactly
$\binom{n_i}{k-\ell}p^{\binom{k-\ell}{2}}$.
Provided we choose $\eta'$ small enough,
this expectation is at least $(1-\eta/2) \tau$, where
$\tau=(np^{\ell})^{k-\ell} p^{\binom{k-\ell}{2}}/(k-\ell)!\sim \mu/M$.
Since $\tau=\Theta(1)$, by a result of Bollob\'as~\cite{BB_Poisson},
the number $N_i'$ of $K_{k-\ell}$s in $G[W_i]$ is asymptotically
Poisson with mean $\tau$. Indeed, $N_i'$
may be coupled with a Poisson distribution $Z$ with mean $(1-\eta)\mu/M$
so that $N_i'\ge Z$ holds with probability $1-o(n^{-\delta})$.
Furthermore, by the first moment method,
with probability $1-o(n^{-\delta})$, the graph $G[W_i]$ does not contain
two $K_{k-\ell}$s sharing a vertex, so $N_i=N_i'$.

It remains only to prove that $\evE_i$ does indeed hold with high
conditional probability.
Recall that at the start of stage $i$, all we know about the
edges between $C_i$ and $V_i$ is that none of the events $\evE_{v,j}$, $v\in
C_i$, $j<i$ holds. This information may be regarded as a separate
condition $\evF_v$ for each $v\in C_i$, where $\evF_v=\bigcap_{j<i}\evE_{v,j}^\cc$
depends only on edges between $v$ and $V_i$.
Given this information, the events $\evE_{v,i}$ are independent,
and each holds with probability
$r=\Pr(\evE_{v,i}\mid \evF_v)$.
Now $\evE_{v,i}$ is an up-set and $\evF_v$ is a down-set, so $r\le
\Pr(\evE_{v,i})=p^{\ell}$. Hence, whatever $|C_i|$ is, the conditional probability
that $n_i\ge 2p^{\ell}n$ is exponentially small. Since
$|C_{i+1}|=|C_i|-n_i$, and we stop after at most $(\log n)^{5k^3}$ steps,
we may thus assume in what follows that $|C_i|\ge n-o(n)$.

Regarding the sets $S_j$, $j\le i$, as fixed, and forgetting our
present conditioning, if all we assume about the edges from
$v$ to $V_i$ is that $\evE_{v,i}$ holds, i.e., that
all edges from $v$ to $S_i$ are present, then each $\evE_{v,j}$, $j<i$,
has conditional probability $p^{|S_j\setminus S_i|}\le p$.
Recalling that we abandon our exploration after at most $(\log n)^{5k^3}$
steps, it follows that
\[
 \Pr(\evF_v^\cc\mid \evE_{v,i}) = \Pr\left(\bigcup_{j<i} \evE_{v,j}\mid \evE_{v,i}\right)
  \le \sum_{j<i}\Pr(\evE_{v,j}\mid \evE_{v,i}) \le ip \le
 (\log n)^{5k^3}p \le \eta',
\]
if $n$ is large enough.
Hence $\Pr(\evF_v\mid \evE_{v,i})\ge 1-\eta'$. In other words, $\Pr(\evF_v\cap \evE_{v,i})\ge
(1-\eta')\Pr(\evE_{v,i})$. This trivially implies that
$\Pr(\evF_v\cap \evE_{v,i})\ge (1-\eta')\Pr(\evF_v)\Pr(\evE_{v,i})$, i.e., that
$\Pr(\evE_{v,i}\mid \evF_v)\ge (1-\eta')\Pr(\evE_{v,i}) = (1-\eta')p^{\ell}$.

It follows that $n_i$ stochastically dominates a Binomial distribution
with parameters $|C_i|$ and $(1-\eta')p^{\ell}$. Since $|C_i|\ge n-o(n)$,
we get the required lower bound on $n_i$, completing the proof.
\end{proof}

Let $N$ denote the number of $K_k$s in $G=G(n,p)$ for which the corresponding
event $\evQ$ holds, and let $N'=N_{\ge (\log n)^{5k^3}}(\Gkl_p)$ be the number
of $K_k$s in {\em large} components of $\Gkl_p$,
that is, components containing at least $(\log n)^{5k^3}$ copies of $K_{\ell}$.
If $V_0$ spans a $K_k$ for which $\evQ$ holds, then by definition
the corresponding component of $\Gkl_p$ contains a tree
with at least $(\log n)^{5k^3}$ vertices; furthermore, exploring
this tree from the root, for each new vertex we find $M=\binom{k}{\ell}-1\ge 1$
new $K_\ell$s. Hence the component is large, so
$N\le N'$.

\begin{lemma}\label{start}
Fix $\eps>0$, and let $p=\Theta(n^{-\frac{2}{k+\ell-1}})$ be chosen 
so that $\mu$ is a constant greater than $1$.
Then
\[
 (\sigma(\mu)-\eps) \nv \le N\le N'
 \le (\sigma(\mu)+\eps) \nv
\]
holds whp.
\end{lemma}
\begin{proof}
Fixing a set $V_0$ of $k$ vertices of $G=G(n,p)$,
recall that $\evA=\evA(V_0)$ is the event that $V_0$ spans a $K_k$
in $G$.
We have $\E(N)=\binom{n}{k}\Pr(\evA)\Pr(\evQ\mid \evA)$, which is
at least $(\sigma(\mu)-o(1))\binom{n}{k}p^{\binom{k}{2}}
=(\sigma(\mu)-o(1))\nv$ by Lemma~\ref{el}.
As noted above, $N\le N'$ always holds.
Thus $\E(N^2)\le \E((N')^2)$.
But $\E((N')^2)\le (\sigma(\mu)^2+o(1))\nv^2$
by Lemma~\ref{ll2}. Hence $\E(N^2)\le (1+o(1))\E(N)^2$,
which implies that $N$ is concentrated around its mean.
Furthermore, $\E(N)\le \E(N')\sim \sigma(\mu)\nv$,
so we have $\E(N)\sim \sigma(\mu)\nv$,
and the result follows.
\end{proof}

\begin{remark}
It is perhaps interesting to note that there is an alternative
proof of the bounds on $N'$ given in Lemma~\ref{start}, using the a sharp-threshold
result of Friedgut~\cite{Friedgut} instead of the second moment method.
Let us briefly outline the argument.
Let $\eva$ be the event that the number $N'=N_{\ge (\log n)^{5k^3}}(\Gkl_p)$
of $K_k$s in large components satisfies
$N'\ge(\sigma(\mu)-\eps)\nv$.
In the light of the expectation bound given by Lemma~\ref{ll1},
it suffices to prove that $\eva$ holds whp.

We view $\eva$ as an event in the probability space $G(n,p)$, in which
case it is clearly increasing and symmetric.
We shall consider $\Pr_{p'}(\eva)$, the probability that $G(n,p')$
has the property $\eva$. When we do so, we keep the definition of $\eva$
fixed, i.e., the definition of $\eva$ refers (via $\mu$ and $\nv$) to $p$, not to $p'$.

Fix $\eta>0$ such that $\sigma(\mu-\eta)>\sigma(\mu)-\eps/4$.
Applying Lemma~\ref{el} with $p'$ reduced by an appropriate
constant factor, we find that
$\E_{p'}(N')\ge \E_{p'}(N)\ge (\sigma(\mu-\eta)-\eps/4)\binom{n}{k}
(p')^{\binom{k}{2}}$,
which is at least $(\sigma(\mu)-3\eps/4)\nv$
if we choose $p'$ correctly.
Since $N'$ is bounded by the total number of $K_k$s, which is very
unlikely to be much larger than its mean $\nv$,
it follows that $\Pr_{p'}(\eva)$ is bounded away from zero.

Since $p/p'$ is a constant larger than $1$, if $\eva$ has a sharp threshold,
we have $\Pr_p(\eva)\to 1$ as required. Otherwise, Theorem 1.2 of
Friedgut~\cite{Friedgut} applies. We conclude that there is a constant $C$
such that $\Pr_p(\eva\mid \eve)\to 1$, where $\eve$ is the event
that a fixed copy of $K_C$ is present in $G=G(n,p)$.
Of course, conditioning on $\eve$ is equivalent to simply adding the edges
of $K_C$ to $G$. Hence, whp, $G(n,p)$ has the property
that after adding a particular copy of $K_C$ to $G$, the event
$\eva$ holds. But the expected number of $K_k$s in $G\cup K_C$ that
share at least $\ell$ vertices with a $K_k$ in $G\cup K_C$
not present in $G$ turns out to be less that $n^{-\delta}\nu$
for some $\delta>0$.
Hence, $G\cup K_C$ contains at most $n^{-\delta/2}\nv$ such $K_k$s
whp. Whenever this holds, removing the edges of $K_C$ from $G$ splits
existing components into at most $n^{-\delta/2}\nv$ new components.
It follows that $G$ has whp at most $n^{-\delta/2}\nv(\log n)^{5k^3}$
fewer $K_k$s in large components that $G\cup K_C$.
Since $G\cup K_C$ has property $\eva$ whp, it follows that $G$
has the same property with a slightly increased $\eps$ whp.
\end{remark}

At this point, we have shown that whp we have the `right' number
of $K_k$s in `large' components; it remains to show that in fact
almost all such $K_k$s are in a single giant component.
In the special case $k=2$, $\ell=1$, i.e., when $\Gkl_p$ is simply
$G(n,p)$, there are many simple ways of showing this, most
of them based on `sprinkling' of one form or another:
following the original approach of Erd\H os and R\'enyi~\cite{ERgiant}
to the study of the giant component of $G(n,p)$, one chooses $p'$ slightly
smaller than $p$, and views $G(n,p)$ as obtained from $G(n,p')$ by
`sprinkling' in a few extra edges. Using independence of the sprinkled
edges from $G(n,p')$, it is easy to show that whp the sprinkled
edges join up almost all large components of $G(n,p')$ into
a single giant component.
Unfortunately, most of these approaches
do not carry over to the present setting;
the essential problem is that, depending on the parameters, $\Gkl_p$
may well have many more vertices than $G(n,p)$. In fact, it may have
many more than $n^2$ vertices.
Approaches such as forming an auxiliary graph on the large components,
joining two if they are connected by sprinkled edges, and then
comparing this graph to $G(n',p')$ for suitable $n'$ and $p'$,
do not seem to work: here $n'$ is much larger than $n$, and there
is not nearly enough independence for such a comparison
to be possible.
For the same reason, we cannot count cuts between largish components,
and estimate the number not joined by sprinkled edges: we may have
many more than $2^{n^2}$ cuts, while the probability that a given
cut is not joined will certainly be at least $2^{-n^2}$.

Fortunately, we can get another version of the sprinkling argument to work:
the key result is the following rather ugly lemma.
In stating this we write $p_0$ for $n^{-\frac{2}{k+\ell-1}}$,
so $\mu(n,p)=\Theta(1)$ is equivalent to $p=\Theta(p_0)$.
We write $\nv_0$ for $\nv(p_0)=\binom{n}{k}p_0^{\binom{k}{2}}$.

\begin{lemma}\label{sprinkle}
Fix constants $\eps>0$ and $A>0$, let $G_0$ be any graph on $[n]$, and
let $C_1,C_2,\ldots,C_r$ list all components of the corresponding
graph $\Gkl_0$ that contain one or more $K_k$s in $G_0$
with property $\evQ$.
Suppose that
\begin{enumerate}
\item\label{a1} between them the $C_i$ contain
at least $2\eps\nv_0$ copies of $K_k$ in $G_0$,
\item\label{a2} no single $C_i$ contains all but $\eps\nv_0$
copies of $K_k$ in $G_0$,
\item\label{a3}
$G_0$ contains at most $A\nv_0$ copies of $K_k$,
\item\label{a4}
for $1\le s<k$ we have
\[
 Z_s \le A n^{2k-s}p_0^{2\binom{k}{2}-\binom{s}{2}},
\]
where $Z_s$ is the number of pairs of $K_k$s in $G_0$ sharing
exactly $s$ vertices, and
\item\label{a5}
no vertex of $G_0$ lies in more than $\nv_0/\sqrt{n}$
copies of $K_k$ in $G_0$.
\end{enumerate}
Set $\gamma=(\log n)^{-4}$, let $G=G(n,\gamma p_0)$ be a random
graph on the same vertex set as $G_0$, and let $\Gkl_1\supset \Gkl_0$ be the graph
$\Gkl$ derived from $G_1=G_0\cup G$. 
Then, for any fixed $i$, the probability that there is some $j$
such that $C_i$ and $C_j$ are contained in a common component
of $\Gkl_1$ is at least $c$, for some constant $c=c(A,\eps)>0$
depending only on $A$ and $\eps$.
\end{lemma}

In other words, roughly speaking, and ignoring all the conditions for a moment,
sprinkling in extra edges with density $\gamma p_0$ is enough to give
any given `large' component of $\Gkl_0$ at least probability $c$
of joining up with another such component, for some $c>0$ that does
not depend on $n$.

We shall prove Lemma~\ref{sprinkle} later; first, we show that
Theorem~\ref{th1} follows.
\begin{proof}[Proof of Theorem~\ref{th1}]
Let $p=p(n)$ be chosen so that $\mu=\mu(n,p)$ is constant and $\mu>1$.
It suffices to show that for any $\eps>0$,
\begin{equation}\label{lwhp}
 C_1(\Gkl_p) \ge N_0=(\sigma(\mu)-2\eps) \nv
 =(\sigma(\mu)-2\eps) \binom{n}{k}p^{\binom{k}{2}}
\end{equation}
holds whp: letting $\eps\to 0$, \eqref{lwhp} implies that
$C_1(\Gkl_p)\ge (\sigma(\mu)-\op(1))\nv$, while
Lemma~\ref{ll1} immediately implies that
$\E(C_1(\Gkl_p)) \le (\sigma(\mu)+o(1))\nv$.
Together, these two statements imply that
$C_1(\Gkl_p)=(\sigma(\mu)+\op(1))\nv$, which is what
the first statement of Theorem~\ref{th1} claims.
For the second, we simply observe that the same argument gives
$C_1(\Gkl_p)+C_2(\Gkl_p)=(\sigma(\mu)+\op(1))\nv$.

To establish \eqref{lwhp}, let us choose $p'<p$ so that
$(\sigma(\mu(p'))-\eps/3)(p'/p)^{\binom{k}{2}}\ge \sigma(\mu)-\eps$.
From continuity of $\sigma$, we can choose such a $p'$ with $p-p'=\Theta(p_0)$.
By Lemma~\ref{start}, applied with $p'$ in place of $p$ and $\eps/3$
in place of $\eps$, whp at least $N_1=N_0+\eps\binom{n}{k}p^{\binom{k}{2}}$
copies of $K_k$
in $G(n,p')$ have property $\evQ$; let
$V_1,\ldots, V_{N_1}$ be (the vertex sets of) $N_1$ such copies.

Let $T=(\log n)^3$, and let $H_1,\ldots,H_T$ be independent
copies of $G(n,\gamma p_0)$ that are also independent of
$G_0=G(n,p')$, with the vertex sets of $H_1,\ldots,H_T$ and of $G_0$ identical.
Set $G_t=G_0\cup\bigcup_{i=1}^t H_i$,
and note that $G_T$ has the distribution of the random graph $G(n,p'')$
for some $p''$. Since $p'+T\gamma p_0\le p$ if $n$ is large enough,
we have $p''\le p$ if $n$ is large enough, so we may couple $G_T$ and
$G(n,p)$ so that the latter contains the former.
Hence, it suffices to prove that whp there is a single
component of $\Gkl(G_T)$ containing at least $N_0$
of the $k$-cliques $V_1,\ldots,V_{N_1}$.

As the reader will have guessed, we shall sprinkle in edges
in $T$ rounds, applying Lemma~\ref{sprinkle}
successively with each pair $(G_{t-1},H_t)$ in place of $(G_0,G)$,
and $\eps'=\eps \nv/\nv_0 = \eps(p/p_0)^{\binom{k}{2}}$ in place of $\eps$.
As noted above, by Lemma~\ref{start}, whp $G_0$ contains
at least $N_1=(\sigma(\mu)-\eps)\nv$ copies of $K_k$
with property $\evQ$. We may assume that $\eps<\sigma(\mu)/3$, in which
case $N_1\ge 2\eps\nv = 2\eps'\nv_0$.
Since
the event that $V_i$ has property $\evQ$ is increasing, and $G_0\subset G_t$
for all $t$, whp the first assumption of Lemma~\ref{sprinkle}
holds for $G_0$ and hence for all $G_t$.

If the second assumption
fails at some point, then we are done: $G_t$ and hence $G_T\supset G_t$
already contains a single component containing at least
$N_1-\eps'\nv_0=N_1-\eps\nv=N_0$ copies of $K_k$, as
required.
The remaining assumptions are down-set conditions, bounding the number
of copies of certain subgraphs in $G_t$ from above.
Standard results tell us that $G(n,p)$ satisfies these
conditions whp if we choose $A$ large enough;
it follows that whp $G_T$ and hence every $G_t$
does too.

From the comments above, we may assume that the conditions of
Lemma~\ref{sprinkle} are satisfied at each stage. Suppose that after
$t$ rounds, i.e., $t$ applications
of Lemma~\ref{sprinkle}, the sets $V_1,\ldots,V_{N_1}$ are now contained
in $r=r(t)$ components $C_1,\ldots,C_r$ of $\Gkl_t$.
By Lemma~\ref{sprinkle}, each $C_i$ has a constant probability $c>0$
of joining
up with some other $C_j$ in each round, so after $(\log n)^2$ further
rounds, the probability that a particular $C_i$ has not joined some
other $C_j$ is at most $(1-c)^{(\log n)^2}=o(n^{-2})$.
It follows that with probability $1-o(n^{-1})$, after $(\log n)^2$
rounds every $C_i$ has joined some other $C_j$.
If this holds, the number $r'$ of components
containing $V_1,\ldots,V_{N_1}$ is now at most $r/2$.
Hence, after $\log r\le \log n$ sets of $(\log n)^2$ rounds,
either an assumption is violated, or there is a single
component containing all $V_i$. But as shown above, there is only one
assumption
that can be violated with probability bounded away from zero, and if
this assumption is violated at some stage, we are already done.
\end{proof}

It remains only to prove Lemma~\ref{sprinkle}.

\begin{proof}[Proof of Lemma~\ref{sprinkle}.]
We assume without loss of generality that $i=1$.
Let $a=\ceil{(\log n)^{5k^3}}$. Since $C_1$ contains a $K_k$ with property
$\evQ$, $C_1$ contains at least $a$ distinct copies of $K_{\ell}$,
each lying in a $K_k$ in $C_1$. Let $S_1,\ldots,S_a$
be $a$ such copies.

From Assumptions~\ref{a1} and~\ref{a2}, $C_2,\ldots,C_r$ between them contain at
least $\eps\nv_0$ copies of $K_k$.
The set $V_0=V(S_1)\cup \cdots\cup V(S_a)$ has size $O(a)=O^\star(1)$,
and so, using Assumption~\ref{a5},
meets at most $o(\nv_0)$ copies of $K_k$ in $G_0$.
It follows that we may find $b=\eps\nv_0/3$
copies $D_1,\ldots,D_b$ of $K_k$ in $C_2,\ldots,C_r$ such
that each $D_j$ is
vertex disjoint from $V_0$. (We round $b$ up to the nearest integer,
but omit this irrelevant distraction from the formulae.)

It suffices to show that with probability bounded away from zero,
there is a path of $K_k$s in $\Gkl_1$ joining some $S_i$ to some $D_j$.
We shall do this using the second moment method. For this, it helps
to count only paths with a simple form.

By a {\em potential $k$-path} we mean a sequence $V_1,\ldots,V_k$ of sets
of $k$ vertices of $G_0$ with the following properties:
$V_1$ contains some $S_i$, all other vertices of $\bigcup_{t=1}^k V_t$
lie outside $V_0$ (and hence outside $S_i$),
$V_k$ coincides with some $D_j$,
and for $2\le t\le k$, $V_t$ consists of $k-\ell$ vertices
outside $\bigcup_{1\le s<t} V_s$ together with $\ell$ vertices of
$V_{t-1}$, not all of which lie in $\bigcup_{s<t-1}V_s$.

A potential $k$-path starting at $S_i$ and ending at $D_j$ contains
exactly $k(k-\ell)-k$ vertices outside $S_i\cup D_j$: starting with $S_i$
we add $k-\ell$ new vertices for each set $V_t$ in the path, but this count
includes the vertices of $D_j$.
It follows that the number of potential $k$-paths joining $S_i$ to $D_j$
is $\Theta(n^{k(k-\ell)-k})$, so the total number of potential $k$-paths
is $\Theta(abn^{k(k-\ell)-k})$.

A potential $k$-path $(V_t)_{t=1}^k$ joining $S_i$ to $D_j$
is a {\em $k$-path} if all edges contained in each $V_t$ but not 
in $S_i$ or $D_j$ are present in $G=G(n,\gamma p_0)$.
Note that any potential $k$-path contains exactly
$r=k(\binom{k}{2}-\binom{\ell}{2})-\binom{k}{2}$ such edges:
for each $t$ there are $\binom{k}{2}-\binom{\ell}{2}$ edges spanned by $V_t$
but by no earlier $V_s$, but this count includes all edges of $D_j$.

Let $X$ denote the number of $k$-paths.
If any $k$-path is present, then some $S_i$ is joined to some $D_j$
in the graph $\Gkl_1$ formed from $G_0\cup G$, so it suffices
to show that $\Pr(X>0)$ is bounded away from zero.

Since each potential $k$-path is present with probability exactly $(\gamma p_0)^r$,
we have
\begin{eqnarray*}
 \E(X) &=& \Theta\left(abn^{k(k-\ell)-k}(\gamma
   p_0)^{k(\binom{k}{2}-\binom{\ell}{2})-\binom{k}{2}}\right) \\
 &=& \Theta\left( ab \left(n^kp_0^{\binom{k}{2}}\right)^{-1} \left(n^{k-\ell} p_0^{\binom{k}{2}-\binom{\ell}{2}}\right)^k\gamma^r\right).
\end{eqnarray*}
Now the bracket raised to the power $k$ in the last
line above is $\Theta(1)$ by definition of
$p_0$, while $b=\eps \nv_0/3 =\Theta\bb{n^kp_0^{\binom{k}{2}}}$. Thus
we have $\E(X)=\Theta( a \gamma^r)$.
Since, crudely, $r\le k\binom{k}{2}\le k^3/2$, while $a\ge (\log n)^{3k^3}$
and $\gamma=(\log n)^{-4}$,
we have $\E(X)\to \infty$.

It remains to estimate the second moment of $X$. For this, it turns out to
be easier to consider a related random variable $Y$.

A {\em potential free $k$-path} is defined exactly as a potential $k$-path,
except that we omit the condition that $V_k$ coincides with some $D_j$.
It is easy to see that the fraction of potential free $k$-paths that
are potential $k$-paths is exactly $b/\binom{n-|V_0|}{k}=\Theta(b/n^k)=
\Theta\bb{p_0^{\binom{k}{2}}}$.

A {\em free $k$-path} is a potential free $k$-path in which all edges
except those contained in the starting set $S_i$ are present in $G=G(n,\gamma
p_0)$. Note that there are $r'=r+\binom{k}{2}$ such edges, so each potential free
$k$-path is an actual free $k$-path with probability $(\gamma p_0)^{r+\binom{k}{2}}$.
Let $Y$ denote the number of free $k$-paths. It follows that
\begin{equation}\label{YX}
 \E(Y) = \Theta\left(\E(X) p_0^{-\binom{k}{2}} (\gamma
   p_0)^{\binom{k}{2}}\right)
 = \Theta\left(\E(X)\gamma^{\binom{k}{2}}\right) = \Theta(a\gamma^{r'}).
\end{equation}

For $0\le s\le k$, let $Z_s$ denote the number of ordered pairs of copies of $K_k$
in $G_0$ sharing exactly $s$ vertices, and let $Z_s'\le Z_s$ denote the number
of such pairs lying entirely outside $V_0$.
Let $X_s$ denote the number of ordered pairs of $k$-paths whose destinations
(final sets $V_k$) share exactly $s$ vertices, and $Y_s$ the number
of ordered pairs of free $k$-paths with this property.
Among ordered pairs $(P_1,P_2)$ of potential free $k$-paths
whose destinations share $s$ vertices, the fraction of pairs
in which $P_1$ and $P_2$ are also
potential $k$-paths
is exactly $Z_s' \left(\binom{n-|V_0|}{k}\binom{k}{s}\binom{n-|V_0|}{k-s}\right)^{-1}
= \Theta( Z_s' n^{-(2k-s)} ).$
Moreover, this statement remains true if we restrict our attention
to pairs $(P_1,P_2)$ with a certain number of common edges.
Indeed, under any sensible assumption on $(P_1,P_2)$, the pair $(V_k,V_k')$
of destinations of a random pair $(P_1,P_2)$ is uniform on all pairs
of $k$-sets in $[n]\setminus V_0$ sharing $s$ vertices.

Given a pair of paths with destinations sharing $s$ vertices,
for both paths to be present as free $k$-paths requires the presence
of $2\binom{k}{2}-\binom{s}{2}$ more edges in $G$ than
required by their presence as $k$-paths. It follows that
\[
 \E(X_s)/\E(Y_s) =
 \Theta\left( Z_s' n^{-(2k-s)} (\gamma p_0)^{-2\binom{k}{2}+\binom{s}{2}}\right).
\]
By Assumption~\ref{a4} we have
$Z_s=O\bb{n^{2k-s}p_0^{2\binom{k}{2}-\binom{s}{2}}}$ for $1\le s\le k-1$.
This also holds for $Z_k$ by Assumption~\ref{a3},
and hence also for $Z_0\le Z_k^2$.
Hence, for $0\le s\le k$,
\[
 \E(X_s)/\E(Y_s) = O\left( \gamma^{-2\binom{k}{2}+\binom{s}{2}}\right) =
O\left(\gamma^{-2\binom{k}{2}}\right).
\]
Since $\E(X^2)=\sum_{s=0}^k \E(X_s)$ and $\E(Y^2)=\sum_{s=0}^k \E(Y_s)$, it
follows immediately that $\E(X^2)/\E(Y^2)=O(\gamma^{-2\binom{k}{2}})$.
We claim that $\E(Y^2)=O(\E(Y)^2)$.
Recalling from \eqref{YX} that $\E(X)/\E(Y)=\Theta(\gamma^{-\binom{k}{2}})$,
it then follows that $\E(X^2)=O(\E(X)^2)$, and hence that $\Pr(X>0)$
is bounded away from zero. 

To evaluate $\E(Y^2)$, we could argue from the fact that
free $k$-paths are balanced in a certain sense, but rather than
make this precise, it turns out to be easier to simply use
our coupling results from Subsection~\ref{ss_ub}.

We may evaluate $Y$, and hence $Y^2$, as follows.
Start with our set $V_0$ of `reached' vertices, namely $V_0=\bigcup_{i=1}^a
V(S_i)$.
Also, mark $S_1,\ldots,S_a$ as untested copies of $K_{\ell}$.
Now explore as in the proofs of Theorem~\ref{th1l} and Lemma~\ref{ub},
except that we only look
for new vertices outside $V_0$; note
that our edge probability is now $\gamma p_0$ rather than $\Theta(p_0)$,
so the corresponding branching process is strongly subcritical.
We stop the exploration after $k$ `rounds', in the terminology of
Lemma~\ref{ub}; of course it may well die earlier.

We consider three cases. Firstly, let $\evA$ be the event that in the exploration
just described, we find no exceptional edges. Since $|V_0|=O^\star(1)$,
and the total size of the relevant branching processes is also $O^\star(1)$
whp, we have $\Pr(\evA^\cc)=O^\star(\gamma p_0)=O(n^{-\delta})$ for some $\delta>0$
depending only on $k$ and $\ell$.
When $\evA$ holds, we obtain a coupling of our exploration with
$a$ independent copies of the branching process $\bp(\la)$,
where $\la=\mu(\gamma p_0)=\Theta(\gamma^{\binom{k}{2}-\binom{\ell}{2}})$.
If $\evA$ holds, the number of $K_k$s reached in the final round is equal
to $N_k/M$, where $N_k$ is the number of particles in generation $k$
of the combined branching process, and we divide
by $M=\binom{k}{\ell}-1$ since we add $M$ copies of $K_{\ell}$ for
each $K_k$ we find.

Now from standard
branching process results, $\E(N_k^2)=\Theta(a^2\la^{2k})
=\Theta(a^2\gamma^{2r'})$, recalling that
$r'=k(\binom{k}{2}-\binom{\ell}{2})$ is the number of edges of $G_0$
in a free $k$-path.
It follows that $\E( Y^2 1_\evA) = O(a^2\gamma^{2r'})$.

We claim that there is a constant
$K$ such that the chance of finding more than $K$ exceptional edges
is $o(n^{-10k^3})$. To see this, first note that the probability
that a Poisson random variable with mean at most $1$ exceeds
$\log n$ is of order $(\log n)^{-\log n}=o(n^{-20k^3})$.
Hence, with probability $1-o(n^{-10k^3})$, the first $k$ generations
of $a+\log n$ copies of $\bp(\la)$ contain at most $(a+\log n)k(\log n)^k
 = O^\star(1)$
particles -- simply crudely bound the number of children
of each particle by $\log n$.
Now arguing as in the proof of Theorem~\ref{th1l}, given that we
have reached $O^\star(1)$ vertices, the chance of finding
an exceptional edge is at most $n^{-\delta}$ for some $\delta>0$.
Hence, the chance of finding $K$ such edges within the first $O^\star(1)$
steps is $O^\star(n^{-\delta K})$ which is $o(n^{-10k^3})$ if we pick
$K$ large enough. But if we find no more than $K$ exceptional
edges within $O^\star(1)$ steps, and the first $k$ generations
of $a+K\le a+\log n$ branching processes have total size $O^\star(1)$,
then (recalling that we stop
after $k$ rounds),
our coupling succeeds, with $a+K$ branching processes as the
upper bound.

Let $\evB$ be the event that we do find more than $K$ exceptional edges,
so $\Pr(\evB)=o(n^{-10k^3})$.
The number of pairs of free $k$-paths present in the complete graph
on $K_n$ is easily seen to be at most $n^{2k^3}$,
so we have $\E(Y^2 1_\evB)\le \Pr(\evB) n^{2k^3}=o(n^{-8k^3}) = o(a^2\gamma^{2r'})$.

Finally, let $\evC=(\evA\cup \evB)^\cc$.
If $\evC$ holds then, as above, with very high probability we have
reached $O^\star(1)$ vertices in our exploration. The picture given
by our exploration may be complicated by the exceptional edges,
but $O^\star(1)$ vertices in any case contain $O^\star(1)$ (pairs of)
free $k$-paths, so we have $\E(Y^2 1_\evC) =O^\star(\Pr(\evC))=O^\star(\Pr(\evA^\cc))=o(1)$.

Putting it all together, $\E(Y^2) = \E(Y^2 1_\evA)+ \E(Y^2 1_\evB)+ \E(Y^2 1_\evC)
= O(a^2\gamma^{2r'})$. From \eqref{YX} we thus have $\E(Y^2)=O(\E(Y)^2)$.
As noted earlier it follows that $\E(X^2)=O(\E(X)^2)$, and
thus that $\Pr(X>0)$ is bounded away from 0, as required.
\end{proof}

\subsection{Far from the critical point}\label{ss_ext}

In the previous subsections we focused on the `approximately critical'
case, where $p$ is chosen so that
the expected number of other $K_k$s adjacent to (i.e.,
sharing at least $\ell$ vertices with) a given $K_k$
is of order $1$. In more standard percolation contexts, one can make this
assumption without loss of generality; using monotonicity it follows
that the fraction of vertices in the largest component tends to $0$
or $1$ outside this range of $p$.

Here we do not have such simple monotonicity, because the
number of vertices of $\Gkl_p$ changes as $p$ varies.
However, it is still easy to deduce results for
values of $p$ outside the range $p=\Theta(p_0)$ from those for $p$ inside
this range.

For $p=o(p_0)$, this is essentially trivial; since the property
of $G$ corresponding to $\Gkl$ containing a component of size
at least $C\log n$ is monotone, Theorem~\ref{th1l}
together with concentration of the number of $K_k$s
trivially implies that the largest component of $\Gkl_p$
contains whp a fraction $o(1)$ of the vertices of $\Gkl_p$,
as long as $\nv=\nv(p)=\binom{n}{k}p^{\binom{k}{2}}$,
the expected number of vertices of $\Gkl_p$, grows
faster than $\log n$. When $\nv$ grows slower than $\log n$
(or indeed than $\sqrt{n}$), by estimating the expected number of cliques
sharing one or more vertices it is very easy to check that whp $\Gkl_p$
contains no edges, and thus no giant component (as long as
$\nv$ does tend to infinity).

To handle the case $p/p_0\to\infty$, we use a slightly different argument.
Let $N$ denote the number of pairs of vertex disjoint cliques
in $G(n,p)$ that lie in the same component of $\Gkl_p$.
Let $p=\Theta(p_0)$. Since the expected number of pairs of cliques
in $G(n,p)$ sharing one or more vertices is $o(\nv^2)$,
Theorem~\ref{th1} shows that $\E_p(N)\ge \bb{\sigma(\mu(p))^2-o(1)}\nv^2$,
considering only pairs in the giant component.
Fix two disjoint sets $V_1$, $V_2$ of $k$ vertices of $G(n,p)$,
and let $\pi_p$ be the probability that $V_1$ and $V_2$ are joined in $\Gkl_p$
given that $V_1$ and $V_2$ are cliques in $G(n,p)$.
Then we have $\E_p(N)=\binom{n}{k}\binom{n-k}{k}p^{2\binom{k}{2}}\pi_p
\sim \nv^2\pi_p$.
Hence, whenever $\mu(p)=\Theta(1)$, we have $\pi_p\ge \sigma(\mu(p))^2-o(1)$.

Now $\pi_p$ is the probability of an increasing event (in the product space
corresponding to the $\binom{n}{2}-2\binom{k}{2}$ possible edges outside
$V_0$, $V_1$), and is hence an increasing function of $p$.
Since $\sigma(\mu)\to 1$ as $\mu\to\infty$, it follows
that $\pi_p\to 1$ if $p/p_0\to \infty$.
Thus, the expected number of unconnected pairs of cliques in $\Gkl_p$
is $o(\nv^2)$ whenever $p/p_0\to\infty$. Since the number of cliques
is concentrated around $\nv$, it follows that whp almost all vertices
of $\Gkl_p$ lie in a single component.

\subsection{Near the critical point}

Der{\'e}nyi, Palla and Vicsek~\cite{DPV} suggest that for $\ell=k-1$,
`at the critical point', i.e., when $p=((k-1)n)^{-1/(k-1)}$,
the largest component in $\Gkl_p$ contains roughly $n$ vertices
of $\Gkl_p$, i.e., roughly $n$ $k$-cliques.
This is based both on computer experiments, and on the heuristic
that at the critical point, the giant component in random graphs
is roughly `treelike'. This latter heuristic seems extremely weak:
there is no reason why a treelike structure in $\Gkl_p$ cannot
contain many more than $n$ $k$-cliques. Indeed, one would not
expect whether or not two $k$-cliques share a single vertex to play much
role in the component structure of $\Gkl_p$.

It would be interesting to know whether the observation of \cite{DPV}
is in fact correct, but there are several problems. Firstly, the
question is not actually that natural: why chose exactly this value
of $p$? In $G(n,p)$, it is natural to take $p=1/(n-1)$ (or $p=1/n$; it turns
out not to matter) as `the' critical probability, since in this
case one has at the beginning a very good approximation
by an exactly critical branching process. However, in general there
is a scaling window within which, for example, the largest
and second largest components are comparable in size. For
$G(n,p)$ the window is $p=n^{-1}+O(n^{-4/3})$; see Bollob\'as~\cite{BBevol}
and {\L}uczak~\cite{Lcrit}; see also the book~\cite{BBRG}.
For other random graph models,
establishing the behaviour of the largest component
in and around the scaling window can be very difficult; see, for example,
Ajtai, Koml\'os and Szemer\'edi~\cite{AKS}, Bollob\'as, Kohayakawa and \L uczak~\cite{BKL},
and
Borgs, Chayes, van der Hofstad, Slade and Spencer~\cite{BCvdHSS_I,BCvdHSS_II,BCvdHSS_III}.

In general, one would expect that inside the scaling window, the largest component
would have size of order $N^{2/3}$, where $N$ is the `volume', which here
would presumably be $\nu=\E(|\Gkl_p|)$. Note that this need not contradict
the experimental results of Der{\'e}nyi, Palla and Vicsek~\cite{DPV}:
it may simply be that their choice of $p$ is (slightly) outside the window.

Unfortunately, due to the dependence in the model, it seems likely to be extremely
difficult to establish results about the scaling window, or about the
behaviour at $p=((k-1)n)^{-1/(k-1)}$. The problem is that there are $o(1)$
errors in the branching process approximation discussed above that appear
right from the beginning. On the one hand, for $\ell=k-1$, as soon as we
find a new $K_k$ sharing $k-1$ vertices with an earlier $K_k$,
there is a probability of order $p$ that a single extra `exceptional'
edge is present forming a $K_{k+1}$, and thus forming extra $K_{k-1}$s
from which we need to explore at the next step. In the other direction,
after even one step of our exploration, we have tested whether any vertex
$v$ not so far reached is joined to all vertices in certain $K_{k-1}$s.
The negative information that $v$ is not so joined reduces the probability
that $v$ is joined to any new $K_{k-1}$
slightly; in fact by a factor of $1-\Theta(p)$ for each $K_{k-1}$ previously
tested. To study the scaling window, or the behaviour
at $p=((k-1)n)^{-1/(k-1)}$ or at $\mu(p)=1$, say, one would presumably need
to understand the net effect of these positive and negative
deviations from the branching process to an accuracy much higher
than the size of each effect. This seems a tall order even for
the first few steps in the branching process, let alone when the component
has grown to size $\Theta(N^{2/3})$ or even $\Theta(n)$.

\section{Variants}

In the rest of the paper we
consider several variants of the clique percolation
problem discussed above. In most cases where we can prove results,
the proofs are minor modifications of those above, so to avoid
trying the reader's patience too far we shall only briefly indicate
the changes.

\subsection{Oriented cliques}

Given $n\ge 2$ and $0\le p\le 1$, let $\Gd(n,p)$ be the random directed
graph on $[n]$ in which each of the $n(n-1)$ possible directed edges
is present with probability $p$, independently of the others.
Thus doubled edges are allowed:
edges $\ore{vw}$ and $\ore{wv}$ may both be present (though
this will turn out to be irrelevant), and the simple graph
underlying $\Gd(n,p)$ has the distribution of $G(n,2p-p^2)$.

Let $\Hd$ be a fixed orientation of $K_k$; for the moment
we shall take $\Hd$ to be $\Kdk$, that is, $K_k$ with a linear order:
$V(\Kdk)=[k]$ and $E(\Kdk)=\{\ore{ij}: 1\le i<j\le k\}$.
Given a directed graph $\Gd$, let $V=V_{\Hd}(\Gd)$ denote the set of all copies of $\Hd$
in $\Gd$. To be totally formal, we may take $V$ to be the set of all
subsets of $\binom{k}{2}$ edges of $\Gd$ that form a graph isomorphic
to $\Hd$. If a given set $S$ of $k$ vertices of $\Gd$
contains double edges, then it may span several copies of $\Hd$,
while if $S$ spans no double edges it spans at most 
one copy of $\Hd$. (For orientations
$\Hd$ with automorphisms, the latter statement
would not be true if we considered
injective homomorphisms from $\Hd$. This is the reason for the somewhat
fussy definition of a `copy' of $\Hd$.)

For $1\le \ell\le k-1$, let $\Gdkl=\Gdkl_{\Hd}$ be the graph formed from $\Gd$ as follows:
let the vertex set of $\Gdkl$ be $V=V_{\Hd}(\Gd)$, and join two vertices
if the corresponding copies of $\Hd$ share at least $\ell$ vertices.
Note that two copies may share $k$ vertices (if double edges are involved);
this will turn out to be irrelevant.
Our aim now is to study the emergence (as $p$ varies)
of a giant component in $\Gdkl_p$,
the graph $\Gdkl$ defined on the copies of $\Hd$ in $\Gd(n,p)$.

\subsubsection{Linearly ordered cliques}

We start by restricting our attention to $\Hd=\Kdk$.
With $\ell=k-1$, the study of this model was proposed by
Palla, Farkas, Pollner, Der{\'e}nyi and Vicsek~\cite{PFPDV}, who
predicted a critical point of $p=(nk(k-1))^{-1/(k-1)}$.
As we shall see, this prediction is correct.

Let us consider the component exploration in $\Gdkl_p$ analogous to
that in $\Gkl_p$ described in Section~\ref{sec_cv}.
The typical case is that we are looking for new $\Kdk$s containing
a given $\Kdl$, say $S$, consisting of $S$ together with $k-\ell$
new vertices.
As before, we expect to find a roughly Poisson number of such new $\Kdk$s,
but now the mean is slightly different: in addition
to choosing a set $N$ of $k-\ell$ new vertices, we must consider
the $k!/\ell!$ linear orders
on $S\cup N$ consistent with the order we
already have on $S$. Given $N$ and such an order,
the probability that this particular $\Kdk$ is present is
then $p^{\binom{k}{2}-\binom{\ell}{2}}$ as before.
As in the undirected case, each new $\Kdk$ we find
typically gives rise to $M=\binom{k}{\ell}-1$ new $\Kdl$s to
explore from in the next step.

Let $\mud=\mud(k,\ell,p)$ be given by
\[
 \mud = \left(\binom{k}{\ell}-1\right)\frac{k!}{\ell!}
 \binom{n}{k-\ell} p^{\binom{k}{2}-\binom{\ell}{2}}.
\]
The proof of Theorem~\ref{th1} goes through
{\em mutatis mutandis} to give the result below. One
can also obtain analogues of the undirected results
for the cases $\mud\to 0$ and $\mud\to\infty$; we
omit these for brevity.

\begin{theorem}\label{th1_d}
Fix $1\le \ell<k$ and let $p=p(n)$ be chosen so that $\mud=\Theta(1)$.
Then, for any $\eps>0$, whp we have
\[
 (\sigma(\mud)-\eps) \nv \le C_1(\Gdkl_p) \le
 (\sigma(\mud)+\eps) \nv,
\]
where $\nu=\binom{n}{k}k!p^{\binom{k}{2}}$ is the expected
number of copies of $\Kdk$ in $\Gd(n,p)$.
\end{theorem}
Note that the function $\sigma$ appearing here is the same function as in Theorem~\ref{th1},
but now evaluated at $\mud$ rather than at $\mu$. In particular, $\sigma(\mud)>0$
if and only if $\mud>1$, and the critical point is given by the solution to $\mud=1$.
In the special case $\ell=k-1$, we have $\mud(k,\ell,p)=(k-1)knp^{k-1}$,
so the critical point is exactly as predicted by
Palla, Farkas, Pollner, Der{\'e}nyi and Vicsek~\cite{PFPDV}.

As the proof really follows that of Theorem~\ref{th1} very closely,
we only briefly describe the differences.
The argument in Subsection~\ref{ss_ub} is essentially unmodified;
it is still true that the first $O(1)$ `exceptional' edges
give rise to the addition of $O(1)$ extra $\Kdl$s, arguing
as before using the total degree of new vertices, rather than in- or out-degree,
say.

For the lower bound, we can argue much of the time using the underlying
undirected graph $G$ rather than $\Gd=\Gd(n,p)$. Indeed, when exploring
from a $\Kdl$ $S_i$, say, we let $W_i$ be the set
of `clean' vertices joined in $G$ to every vertex of $S_i$.
We then look for undirected $k$-cliques in $G[W_i]$.
Arguing as before, the number we find can be coupled
to agree (up to a negligible error term) with a Poisson distribution
with the appropriate mean, now
$(1-\eta)\binom{n}{k-\ell}(2p-p^2)^{\binom{k}{2}-\binom{\ell}{2}}$.
Moreover, as before, we may assume that the $k$-cliques we find
are vertex disjoint. Only at this point do we check the orientations
of the $\binom{k}{2}-\binom{\ell}{2}$ new edges involved
in each $k$-clique; the probability that we find one of the $k!/\ell!$
orientations that gives a $\Kdk$ extending $S_i$ is
$(k!/\ell!)(1/2)^{\binom{k}{2}-\binom{\ell}{2}}+o(1)$,
so the number of such $\Kdk$s that we do find may be closely
coupled to a Poisson distribution with mean $\mud$ as required.

Finally, the argument joining up large components goes through
with only trivial modifications to the definitions.

\subsubsection{Cliques with arbitrary orientations}

We now turn out attention to the phase transition in the graph $\Gdkl_p$
defined on the copies of $\Hd$ in $\Gd(n,p)$, where $\Hd$
is some non-transitive orientation of $K_k$.
Perhaps surprisingly, it turns out that something genuinely new
happens in this case.

\begin{figure}[ht]
 \centering
 \input{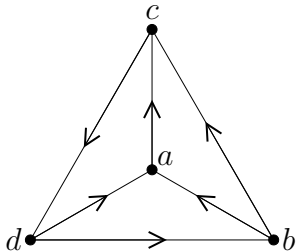}
 \caption{An orientation $\Hd$ of $K_4$.}\label{fig1}
\end{figure}

Let $k=4$, $\ell=3$, let $\Hd$ be the orientation of $K_4$
shown in Figure~\ref{fig1}, and let $\Gdkl_p$ be defined as before.
When exploring a component of $\Gdkl_p$, suppose that we have
found a certain copy of $\Hd$, and are looking for new copies
containing a particular subgraph $S$ of order $3$.
There are now four separate cases, although one
can combine them in pairs.
First suppose the vertex set of $S$ is $\{b,c,d\}$, 
so $S$ is an oriented triangle.
If we find a vertex $v$ joined to $b$, $c$ and $d$,
there are six combinations of orientations
of $vb$, $vc$ and $vd$ that lead to a copy of $\Hd$:
either two edges are oriented towards $v$ and one away, in which
case $v$ plays the role of $a$ in the new copy of $\Hd$,
or two are oriented away from $v$ and one towards,
in which case $v$ plays the role of $b$.
The same holds if $V(S)=\{a,c,d\}$, since $S$ is again
an oriented triangle.

In the other two cases, $S$ is a linearly ordered triangle,
and either $v$ sends edges to the top two vertices
of $S$ and receives an edge from the bottom one, and so plays the role
of $d$, or $v$ sends an edge to the top vertex and receives edges from
the bottom two, playing the role of $c$.

Suppose more generally that $\Hd$ is an orientation of $K_k$
in which no two vertices are equivalent (the orientation
in Figure~\ref{fig1} has this property).
Let $M(\Hd)$ be the $k$-by-$k$ matrix whose $ij$th entry
is the number of ways of orienting the edges
from a new vertex $v$ to $[k]\setminus\{j\}$ such
that $\Hd-j\cup \{v\}$ forms a graph isomorphic to $\Hd$
with $v$ playing the role of vertex $i$.
For example, with $\Hd$ as in Figure~\ref{fig1}, numbering
the vertices in the order $a$, $b$, $c$, $d$, we have
\begin{equation}\label{matrix}
 M = \left(\begin{matrix}
 3 & 3 &  0 & 0 \\
 3 & 3 &  0 & 0 \\
 0 & 0 &  1 & 1 \\
 0 & 0 &  1 & 1 
\end{matrix}\right).
\end{equation}
Let us say that a copy in $\Gd$ of a
subgraph of $\Hd$ induced by $k-1$ vertices
is of {\em type $j$} if it is formed by omitting the vertex $j$.
Also, let us say that a copy of $\Hd$ found in our exploration
by adding a new vertex $v$ to a subgraph of $\Hd$ with $k-1$
vertices is of {\em type $i$} if the new vertex corresponds
to vertex $i$ of $\Hd$.
Then, towards the start of our exploration,
the expected number of type $i$ copies of $\Hd$
we reach from a type $j$ subgraph is $M_{ij}n p^{k-1}$.
When we continue the exploration,
each type $i$ copy of $\Hd$ gives rises to one new subgraph
of each type other than $i$, and this gives us
our branching process approximation.

For the formal statement, less us pass to the general case $1\le \ell\le k-1$.
For simplicity, the reader may prefer to consider only graphs $\Hd$
such that
all $\binom{k}{\ell}$ sets of $k-\ell$ vertices of $\Hd$
are non-equivalent, so that when we extend an $\ell$ vertex graph
to a graph isomorphic to $\Hd$ we can identify which $k-\ell$
vertices of $[k]$ the new vertices correspond to. In general,
we may resolve ambiguous cases arbitrarily. (One could instead collapse
the corresponding types in the branching process, but this
complicates the description.)
Let $M$ be the $\binom{k}{\ell}$-by-$\binom{k}{\ell}$ matrix defined
as follows:
given two $\ell$-element subsets $A$ and $B$ of $[k]$,
let $S$ be the subgraph of $\Hd$ induced by the vertices
in $A$, and consider a set $N$ of $k-\ell$ `new' vertices joined to each
other and to all vertices in $A$. Let $M_{BA}$ be the number of ways
of orienting these new edges so that $A\cup N$ forms a copy of $\Hd$,
and the new vertices correspond to $[k]\setminus B$.
For $\ell=k-1$, this generalizes the definition above.

Let $\bp=\bp_{\Hd}$ be the multi-type Galton--Watson branching process in which
each particle has a type from $\binom{[k]}{\ell}$, started with one particle
of each type, in which children of a particle of type $A$
are generated as follows: first generate independent Poisson
random variables $Z_B$, $B\in \binom{[k]}{\ell}$,
with $\E(Z_B)=M_{BA} \binom{n}{k-\ell} p^{\binom{k}{2}-\binom{\ell}{2}}$.
Then generate $\sum_{B\ne A'} Z_B$ children of each type $A'$.
Let $\sigma_{\Hd}=\sigma_{\Hd}(p)$ denote the survival probability of $\bp$.
The proof of Theorem~\ref{th1_d} extends very easily to prove
the following result.

\begin{theorem}\label{th1_mt}
Fix $1\le \ell<k$ and an orientation $\Hd$ of $K_k$,
and let $p=p(n)$ be chosen so that $n^{k-\ell}p^{\binom{k}{2}-\binom{\ell}{2}}=\Theta(1)$.
Then, for any $\eps>0$, whp we have
\[
 (\sigma_{\Hd}(\mu)-\eps) \nv \le C_1(\Gdkl_p) \le
 (\sigma_{\Hd}(\mu)+\eps) \nv,
\]
where $\nu=\binom{n}{k}(k!/\aut(\Hd))p^{\binom{k}{2}}$ is the expected
number of copies of $\Hd$ in $\Gd(n,p)$.
\noproof
\end{theorem}

Theorem~\ref{th1_mt} is rather unwieldy, but it is not too hard to extract the
critical point. Indeed, in $\bp$ the expected number of type-$B$ children
of a particle of type $A$ is $X_{BA}\binom{n}{k-\ell}p^{\binom{k}{2}-\binom{\ell}{2}}$,
where $X_{BA}=(J-I)M$, with $I$ the identity matrix and $J$ the matrix with all entries
$1$. From elementary branching process results,
the critical value of $p$ is thus given by the solution to
\[
 \lambda \binom{n}{k-\ell}p^{\binom{k}{2}-\binom{\ell}{2}} =1,
\]
where $\lambda$ is the maximum eigenvalue of $X=(X_{BA})$.

Note that this is consistent with Theorem~\ref{th1_d}: taking $\Hd$ to be
$\Kdk$, that is, $K_k$ with a transitive order, it is easy to check
that $M_{BA}=(k-\ell)!$
for every $A,B\in \binom{[k]}{\ell}$. Indeed, we must choose one of
the $(k-\ell)!$ possible orders on the new vertices. Then the relative
order of the new and old vertices is determined by the fact that the new
vertices should play the role of $[k]\setminus B$ in the resulting $\Kdk$.
It follows that $X$ is the $\binom{k}{\ell}$-by-$\binom{k}{\ell}$ matrix
with all entries equal to $(\binom{k}{\ell}-1)(k-\ell)!$,
so
\[
 \lambda=\binom{k}{\ell}\left(\binom{k}{\ell}-1\right)(k-\ell)!=\left(\binom{k}{\ell}-1\right)\frac{k!}{\ell}.
\]

To give a non-trivial application of Theorem~\ref{th1_mt}, let $\Hd$ be the orientation
of $K_4$ shown in Figure~\ref{fig1}. 
Then $M$ is given by \eqref{matrix}, so we have
\[
 X = 
\left(\begin{matrix}
 0 & 1 &  1 & 1 \\
 1 & 0 &  1 & 1 \\
 1 & 1 &  0 & 1 \\
 1 & 1 &  1 & 0 
\end{matrix}\right)
\left(\begin{matrix}
 3 & 3 &  0 & 0 \\
 3 & 3 &  0 & 0 \\
 0 & 0 &  1 & 1 \\
 0 & 0 &  1 & 1 
\end{matrix}\right)
=
\left(\begin{matrix}
 3 & 3 &  2 & 2 \\
 3 & 3 &  2 & 2 \\
 6 & 6 &  1 & 1 \\
 6 & 6 &  1 & 1 
\end{matrix}\right)
\]
It follows that $\lambda$, which may be found as twice
the maximum eigenvalue of a 2-by-2 matrix, is equal to $2(2+\sqrt{13})$,
so the critical $p$ is
$\bb{(4+2\sqrt{13})n}^{-1/3}$.

\subsection{Cliques joined by edges}

In this subsection we return to unoriented graphs,
and consider another natural notion of adjacency for
copies of $K_k$ in a graph $G$: given a parameter $1\le \ell \le k^2$,
two $K_k$s are considered adjacent if they are vertex disjoint and there
are at least $\ell$ edges of $G$ from one to the other.
(One could omit the disjointness condition; much of the time this
will make little difference. Insisting on this condition simplifies
the picture slightly.)
Let $\Gokl(G)$ be the corresponding graph on the copies of $K_k$ in
$G$, and let $\Gokl_p=\Gokl(G(n,p))$ be the graph obtained in
this way from $G(n,p)$. For this notion of adjacency, the most natural special
case to consider is $\ell=1$; the other extreme case, $\ell=k^2$,
of course corresponds to considering copies of $K_{2k}$ sharing $k$ vertices.

It turns out that we can fairly easily determine the percolation
threshold in $\Gokl_p$ for those parameters $(k,\ell)$ for which,
near the threshold, there are `not too many' copies of $K_k$ in $G(n,p)$;
more precisely, there are $o(n)$ copies. This always includes the case
$\ell=1$.

Let $\mu'=\mu'(n,k,\ell,p)$ be given by
\begin{equation}\label{nmudef}
 \mu' = \binom{n}{k}\binom{k^2}{\ell} p^{\binom{k}{2}+\ell},
\end{equation}
and, as before, let $\nv=\nv(n,k,p)=\binom{n}{k}p^{\binom{k}{2}}$
be the expected number of copies of $K_k$ in $G(n,p)$, so
$\nv=\E |\Gokl_p|$.
Let $\bp_0(\la)$ denote a Galton--Watson branching process
in which the offspring distribution is Poisson with mean $\la$,
started with a single particle, and let $\sigma_0(\la)$ denote
the survival probability of $\bp_0(\la)$.
Note that $\sigma_0(\la)n$ is the asymptotic
size (number of vertices) in the largest component of $G(n,\la/n)$.

The following result is analogous to Theorem~\ref{th1},
but, in part due to the extra assumption on $\nv$, much simpler.
\begin{theorem}\label{th_e}
Fix $k\ge 3$ and $1 \le\ell\le k-2$.
Let $p=p(n)$ be chosen so that $\mu'=\Theta(1)$ and $\nv=o(n)$.
Then, for any $\eps>0$,
\begin{equation}\label{nnn}
 (\sigma_0(\mu')-\eps) \nv \le C_1(\Gokl_p) \le
 (\sigma_0(\mu')+\eps) \nv
\end{equation}
holds whp.
\end{theorem}
Note that there is a choice of $p=p(n)$ satisfying
the conditions of Theorem~\ref{th_e} if and only if
$\ell<k/2$. Indeed, the main force of Theorem~\ref{th_e} is to
establish that in this case, the threshold for percolation in $\Gokl_p$
is at the solution $p_0$ to $\mu'(p)=1$, which
satisfies $p_0=\Theta(n^{-\frac{2k}{k(k-1)+2\ell}})$,
with the constant given by \eqref{nmudef}.
As in Section~\ref{sec_cv}, the proof of Theorem~\ref{th_e} will give
an $O(\log n)$ bound in the subcritical case,
as well as an $o(\nv)$ bound on the 2nd largest component
in the supercritical case. The former applies
also for $\ell\ge k/2$, but only under the assumption that $\nv=o(n)$,
i.e., well below what is presumably the critical point in this case.
One can also extrapolate to the highly supercritical case
as in Subsection~\ref{ss_ext}. Here one needs
the condition $\nv=o(n)$ only for the starting value of $p$,
and the conclusion is that for $1\le \ell<k/2$
and any $p$ with $p/p_0\to\infty$ one has, as expected,
almost all vertices of $\Gokl_p$ in a single component.

After these remarks, we turn to the proof of Theorem~\ref{th_e}.
\begin{proof}
We start with the upper bound.
Let us call a copy of $K_k$ in $G(n,p)$ {\em isolated} if it shares
no vertices with any other copies of $K_k$.
Let $N$ and $M$ denote the number of isolated and non-isolated
copies of $K_k$ in $G(n,p)$.
By a standard calculation, the probability that a given copy of $K_k$ is not isolated is
$(1+o(1))k\nu/n=o(1)$, so $\E(M)=o(\nu)$, and whp we have $M=o(\nu)$.
More precisely, we may choose some $\omega=\omega(n)\to\infty$
so that the event $\evB$ that $M\ge \nv/\omega$ has probability $o(1)$.
Since the number of copies of $K_k$ in $G(n,p)$ is concentrated
about its mean, choosing $\omega$ suitably, the event $\evA$ that $|N-\nv|\le
\nu/\omega$ also holds whp.

Let $S_1,S_2,\ldots,S_N$ list the vertex sets of all isolated copies of $K_k$
in $G(n,p)$, and $T_1,\ldots,T_M$ those of all non-isolated copies.
We condition on $N$, $M$, and the sequences $(S_i)$ and $(T_i)$.
We assume that $\evA\setminus \evB$ holds; we may do so since $\Pr(\evA\setminus \evB)=1-o(1)$.
Let $\evE$ denote one of the specific events we condition on,
and let $E^+$ denote the set of all edges lying 
within some $S_i$ or $T_i$, and $E^-$ the set of
all $\binom{n}{2}-|E^+|$ remaining potential edges of $G$.
Let us call a non-empty set $F\subset E^-$ {\em forbidden} if by adding
zero of more edges of $E^+$ to $F$ one can form a $K_k$;
we write $\F$ for the collection of forbidden sets.
The event $\evE$ may be represented as the intersection
of an up-set condition $\evU$, that every edge in $E^+$ is present
in $G(n,p)$, and a down-set condition $\evD$, that no forbidden set is present in $E^-$.
Note that $\evD$ may be regarded as a down-set in $\{0,1\}^{E^-}$.

For the moment, we condition only on $\evU$.
To be pedantic (while, at the same time, committing the common
abuse of using the same notation for a random variable and its possible
values), we fix sequences $(S_i)$ and $(T_i)$ consistent
with $\evA\setminus \evB$, and condition on the event $\evU=\evU((S_i),(T_i))$.
Since we are conditioning only on the presence of a fixed set
of edges, every edge of $E^-$
is present independently with probability $p$.
Let $H$ be the auxiliary graph with vertex set $[N]$ in which
$i$ and $j$ are joined if $S_i$ and $S_j$ are joined by at least
$\ell$ edges.
The probability $p'$ of this event satisfies
\[ 
 p'=\binom{k^2}{\ell}p^{\ell} + O(p^{\ell+1}) \sim \mu'/\nv \sim \mu'/N.
\]
Since, given $\evU$, $H$ has exactly the distribution of $G(N,p')$,
it follows from the classical result of Erd\H os and R\'enyi~\cite{ERgiant}
that whp the largest component of $H$ has order within
$\eps N/2$ of $\sigma_0(\mu')N$.
Note that this corresponds to the desired number of $K_k$s in the
largest component $C$ of $\Gokl_p$. The problem is that we have
not yet conditioned on $\evD$, or allowed for the possible presence
of non-isolated $K_k$s in $C$.

To prove the upper bound in \eqref{nnn} we must account for the
non-isolated $K_k$s. Let us say that $S_i$ and $T_j$ form a {\em bad pair}
if they are joined by $\ell$ edges in $G(n,p)$. Given $\evU$,
the probability of this event is exactly $p'$,
so the expected number of bad pairs $(S_i,T_j)$ is $p'NM=o(p' N^2)=o(N)$.
Similarly, $T_i$ and $T_j$ form a bad pair if they are vertex
disjoint, and joined by at least $\ell$ edges.
The expected number of bad pairs $(T_i,T_j)$ is at most $p'M^2=o(N)$.
Let $H'\supset H$ be the graph on $[N+M]$ defined in the natural
way: two vertices are joined if the corresponding
copies of $K_k$ are disjoint and joined by at least $\ell$ edges.
We have shown that $E(H')\setminus E(H)$, which is exactly
the number of bad pairs, has expectation $o(N)$.

It is well known that for $\la$ fixed, the giant component in $G(m,\la/m)$
is {\em stable upwards}, in the sense that adding $\op(m)$ vertices and
edges cannot increase its size by more than $\op(m)$.
Indeed, this follows from the qualitative form of the distribution
of the small components: for details, see, for example, Theorem 3.9
of Bollob\'as, Janson and Riordan~\cite{BJR}, where the corresponding result is proved
for a more general model. (This result also shows `downwards stability', which
we do not need here. Downwards stability is much harder to prove:
Luczak and McDiarmid~\cite{LMcD} established this for the Erd\H os--R\'enyi model;
in~\cite{BJR}, their argument is extended to the more general model considered there.)
Applying this stability result to $H$, we deduce
that, given $\evU$, we have $C_1(H')=C_1(H)+\op(N)$.
For any $n'$, we have
\[
 \Pr(C_1(H')\ge n' \mid \evE)  = \Pr(C_1(H')\ge n'\mid \evU\cap \evD)
 \le \Pr(C_1(H')\ge n'\mid \evU),
\]
where the inequality is from Harris's Lemma applied in $\{0,1\}^{E^-}$.
Since $H'$ is exactly the graph $\Gokl_p$, the upper bound in \eqref{nnn}
follows.

\smallskip
Turning to the lower bound, we may now ignore the complications
due to non-isolated $K_k$s, and confine our attention to $H$.
However, we must now show that conditioning on $\evD$,
which tends to decrease $C_1(H)$, does not do so too much.
We shall use the same type of argument as in the proof of Lemma~\ref{el}:
exploring $H$ step by step, we shall show that conditioning
on $\evD$ does not decrease the probability of finding an edge in $H$
by showing that finding an edge in $H$
would not decrease the probability of $\evD$ much. There
will be some complications due, for example,
to the possible presence of $K_k$s
made up of edges in $G=G(n,p)$ corresponding to edges in $H$.

As before, we shall condition on $(S_i)$ and $(T_i)$,
assuming that $\evA\setminus \evB$ holds, i.e., that $N\sim \nv$ and
$M=o(\nv)$. In fact, we shall impose a further condition.
Let $\evB'$ be the event that there is a vertex of $G(n,p)$ in more
than $(\log n)^2$ copies of $K_k$, noting
that whether or not 
$\evB'$ holds is determined by the sequences $(S_i)$ and $(T_i)$.
Since $\nu=\binom{n}{k}p^{\binom{k}{2}}=o(n)$,
it is easy to check that $\Pr(\evB')=o(1)$: we omit
the standard calculation which is based on the fact that $K_k$ is strictly
balanced, so having found a moderate number of $K_k$s containing
a given vertex $v$ does not significantly increase the chance of finding
a further such $K_k$.

From now on we condition on the sequences $(S_i)$ and $(T_i)$, 
assuming as we may that $\evA\setminus (\evB\cup \evB')$ holds.
Defining $\evU=\evU((S_i),(T_i))$ and $\evD=\evD((S_i),(T_i))$ as before,
this is again equivalent to conditioning on $\evU\cap \evD$. As before,
since we fix $(S_i)$ and $(T_i)$, the event $\evU$ is simply the event
that every edge in the fixed set $E^+=\bigcup E(S_i)\cup \bigcup E(T_i)$
is present in $G(n,p)$. Note for later that, since $\evB'$ does not hold, we have
\begin{equation}\label{dE+}
 d_{E^+}(v)\le k(\log n)^2
\end{equation}
for every $v\in V(G)$, where $d_{E^+}(v)$ is the number of edges
of $E^+$ incident with $v$.

Let $f_1,f_2,\ldots$ be the $\binom{N}{2}$ possible edges of $H$, listed in an arbitrary
order.
We now describe an algorithm that reveals a subgraph $H_0$ of $H$.
During step $r$, $1\le r\le \binom{N}{2}$,
we shall test whether $f_r$ is present
in $H$, except that if $f_r$, together with some previously
discovered edges of $H_0$, would form a cycle in $H_0$, or
would cause the degree of some vertex of $H_0$ to exceed
$(\log n)^2$, then we omit step $r$.
Step $r$ consists of a series of sub-steps: in each we consider
one of the $\binom{k^2}{\ell}$ sets $I$
of $\ell$ potential edges of $G=G(n,p)$ whose presence would give rise
to the edge $f_r$ in $H$, and test whether all edges in $I$ are present
in $G$. If such a test succeeds, we add $f_r$ to $H_0$,
and omit further tests for the
same $f_r$, i.e., continue to step $r+1$.

Suppose that we have reached the $t$th sub-step of the algorithm described
above, and let $I=I_t$ be the set of $\ell$ potential edges of $G$
whose presence we are about to test for.
We claim that, given the history, the conditional probability
that all edges in $I$ are is present is $(1+o(1))p^\ell$.
More precisely, let $E_t^+$ be the union of all sets $I_s$, $s<t$,
which we found to be present, and let $\evU_t=\{E_t^+\subset E(G)\}$.
Also, let $\F_t$ be the set
of sets $I_s$, $s<t$, found to be absent,
and let $\evD_t$ be the event that no $F\in \F_t$ is present in $E(G)$.
Recalling that we start by conditioning on $\evU\cap \evD$,
the algorithm reaches its particular present state if and only if
$\evU\cap \evD\cap \evU_t\cap \evD_t$ holds, so our precise claim is that
for any $\eta>0$,
if $n$ is large enough, then for any possible $I_t$, $\evU_t$ and
$\evD_t$ we have
\begin{equation}\label{cl}
 \Pr\bb{I_t\subset E(G) \mid \evU\cap \evD\cap \evU_t\cap \evD_t}\ge (1-\eta)p^\ell.
\end{equation}

Before proving \eqref{cl}, let us see that the Theorem~\ref{th_e} follows.
Let $H_1$ be the union of $H_0$ and all edges $f_r$ which we
omitted to test. Assuming \eqref{cl}, we always have
\begin{multline}\label{dom}
 \Pr\bb{f_r\in E(H_1) \mid E(H_1)\cap\{f_1,\ldots,f_{r-1}\} } \\
 \ge (1-\eta-o(1))\binom{k^2}{\ell}p^\ell \sim (1-\eta)p'\sim (1-\eta)\mu'/N.
\end{multline}
Indeed, if $f_r$ is omitted, the conditional probability above
is $1$ by definition; otherwise, we apply \eqref{cl} to the
$\binom{k^2}{\ell}$ sub-steps associated to $f_r$.
Now \eqref{dom} tells us that for $n$ large enough,
$H_1$ stochastically dominates $G(N,(1-2\eta)\mu'/N)$, say.
Taking $\eta$ small enough, it follows that whp
\begin{equation}\label{CH1}
 C_1(H_1)/N\ge \sigma_0((1-2\eta)\mu')-\eps/4 \ge \sigma_0(\mu')-\eps/2.
\end{equation}

If $\Delta(H)\le (\log n)^2-1$, then we only omit step $r$
if adding $f_r$ would create a cycle, so in this case
$H_0$ is the union of one spanning tree for each component
of $H$, and all edges of $H_1$ join vertices of
$H_0$ that are already joined by paths in $H_0$.
Hence $C_1(H_0)=C_1(H)=C_1(H_1)$. As noted earlier,
given only $\evU$, the graph $H$ has exactly the distribution
of $G(N,p')$. Since $\evU$ is a principal up-set, and $\evD$ is a down-set,
it follows that the distribution
of $H$ given $\evU\cap \evD$, which is what
we are considering here, is stochastically dominated by that of
$G(N,p')$. Since $Np'\sim \mu'=\Theta(1)$, it follows
that whp $\Delta(H)\le (\log n)^2-1$, so whp $C_1(H)=C_1(H_1)$.
Since $N\sim \nu$, this together with \eqref{CH1} gives
the lower bound in \eqref{nnn}, completing the proof
of Theorem~\ref{th_e}.

\smallskip
It remains only to prove \eqref{cl}. Let us start by observing that
$E_t^+\cup I_t$ cannot contain any forbidden
set $F\in \F$, i.e., that the set $E^+\cup E_t^+\cup I_t$ 
contains no $K_k$ other than $S_1,\ldots,S_N, T_1,\ldots,T_M$. 
This is true of $E^+\cup E_t^+$ since we condition on $\evD$,
and we are assuming that the present state of the algorithm
is a possible one. Suppose then that adding $I_t$ to $E^+\cup E_t^+$
creates a new copy $K$ of $K_k$, and let the edge
$f_r$ we are testing be $ij$. 
Now $E_t^+$ contains
edges between $S_{i'}$ and $S_{j'}$ if and only if we have already
found the edge $i'j'$ in $H_0$. If $K$ meets
three or more of the $S_{i'}$, then $H_0\cup f_r$ would contain
a triangle, which is impossible by definition of the algorithm.
This leaves only the case that $K$ meets exactly two sets $S_{i'}$,
which must be $S_i$ and $S_j$. But then the only edges of $E_t^+\cup I_t$
between $S_i$ and $S_j$ are those of $I_t$.
Now $K$ contains at least $k-1$ edges between these sets,
while $|I_t|=\ell<k-1$, so there is no such $K_k$.

There are two types of conditioning in \eqref{cl}, that on $\evU\cap \evU_t$
and that on $\evD\cap \evD_t$. The first type is trivial, since
$\evU\cap \evU_t$ is simply the event that every edge in $E^+\cup E_t^+$
is present.
Let $X=E(K_n)\setminus (E^+\cup E_t^+)$. Then we may as well
work in $\Pr^X$, the product probability measure on $\{0,1\}^X$
in which each edge is present with probability $p$.
Let 
\begin{equation}\label{tFdef} 
\tF_t=\{F\cap X: F\in \F\cup \F_t\},
\end{equation}
and let $\tD\subset \{0,1\}^X$
be the event that none of the `forbidden' sets in $\tF_t$ is present,
so $\Pr^X(\tD)=\Pr(\evD\cap \evD_t\mid \evU\cap \evU_t)$.
Also, let $\evIt\subset \{0,1\}^X$ be the event that all edges in $I_t$
are present, noting that $I_t\subset X$.
Then \eqref{cl} is equivalent to
\begin{equation}\label{cl2}
 \Pr^X(\evIt\mid \tD_t)\ge (1-\eta)\Pr^X(\evIt).
\end{equation}
The key idea is to split $\tF_t$ into two sets, $\tF'$ and $\tF''$,
the first consisting of those $F$ that intersect $I_t$,
and the second those that do not. It turns out that we can ignore
the ones that do not.
More precisely, let $\tD'$ be the event that no $F\in \tF'$ is present,
and $\tD''$ the event that no $F\in \tF''$ is present, so $\tD_t=\tD'\cap \tD''$.

We may rewrite \eqref{cl2} in any of the following forms, which
are step-by-step trivially equivalent: (we drop the superscript $X$
at this point, since the events we are now considering are in any
case independent of edges outside $X$)

\begin{eqnarray}
 \Pr(\evIt\mid \tD'\cap \tD'') &\ge& (1-\eta) \Pr(\evIt)    \nonumber\\
 \Pr(\evIt\cap \tD'\cap \tD'') &\ge& (1-\eta) \Pr(\evIt)\Pr(\tD'\cap \tD'') \nonumber\\
 \Pr(\tD'\cap \tD''\mid \evIt) &\ge& (1-\eta) \Pr(\tD'\cap \tD'') \nonumber\\
 \Pr(\tD'\mid \evIt\cap \tD'')\Pr(\tD''\mid \evIt) &\ge& (1-\eta) \Pr(\tD'\mid \tD'')\Pr(\tD'') \nonumber\\
 \Pr(\tD' \mid \evIt\cap \tD'') &\ge& (1-\eta)\Pr(\tD'\mid \tD''). \label{last}
\end{eqnarray}
The only step which is not trivial from the definition of conditional
probability is the last one: for this we note that by definition $\tD''$
depends only on edges of $X\setminus I_t$.

We shall prove \eqref{last} by simply ignoring the conditional
probability on the right,
showing that
\begin{equation}\nonumber% \label{last2}
 \Pr(\tD'\mid \evIt\cap \tD'') \ge (1-\eta).
\end{equation}
This clearly implies \eqref{last}, and hence, from the argument above,
implies \eqref{cl}.
Let $\tU'$ be the complement of $\tD'$, so our aim is to show
that $\Pr(\tU'\mid \evIt\cap \tD'')\le \eta$.
Now $\evIt$ is again an event of a very simple form, that a certain particular
set $I_t$ of edges is present.
Since $\tU'$ is an up-set while $\tD''$ is a down-set,
applying Harris's Lemma in $\{0,1\}^{X\setminus I_t}$, it follows that
\[
 \Pr(\tU'\mid \evIt\cap \tD'')  \le \Pr(\tU'\mid \evIt),
\]
so it suffices to prove that $\Pr(\tU'\mid \evIt)\le \eta$.

At this point we have eliminated all non-trivial conditioning;
all that is left is counting. Indeed,
\begin{equation}\label{sss}
 \Pr(\tU'\mid \evIt) \le \sum_{F'\in \tF'} p^{|F'\setminus I_t|}.
\end{equation}
Recalling \eqref{tFdef}, there are two contributions to the sum above.
The first is from sets $F'\in \tF'$ corresponding to sets $I_s\in \F_t$,
i.e., to failed tests
for previous $I_s$. By definition of $\tF'$, we have such
an $F'\in \tF'$ if and only if $I_s\cap I_t\ne \emptyset$,
in which case $I_s$ and $I_t$ correspond to the same potential edge $f_r$ of $H$.
But then there are at most $\binom{k^2}{\ell}-1$ possibilities
for $I_s$, and for each we have $|F'\setminus I_t|=|I_s\setminus I_t|\ge 1$,
so the contribution to \eqref{sss} is $O(p)=o(1)$.

The remaining terms come from $F'=F\cap X$ with $F\in \F$
and with $F'\cap I_t\ne\emptyset$, i.e., with $F\cap I_t\ne \emptyset$.
Recalling that $F$ is a set of edges that, together with $E^+$,
would create a $K_k$, it thus suffices to show that
\begin{equation}\label{last3}
 \sum_K p^{E(K)\setminus (E^+\cup E_t^+\cup I_t)} =o(1),
\end{equation}
where the sum runs over all copies of $K_k$ on $V(G)$ containing
at least one edge from $I_t$. 
Now $H_0$ has maximum degree at most $(\log n)^2$ by the definition
of our algorithm. Hence
$d_{E_t^+}(v)\le \ell(\log n)^2$ for every $v\in V(G)$.
Using  \eqref{dE+} it follows that
the graph $G'$ on $V(G)$ formed by the edges in $E^+\cup E_t^+\cup I_t$
has maximum degree at most $(k+\ell)(\log n)^2+\ell \le (\log n)^3$, say.
This is all we shall need to prove \eqref{last3}.

Let $Z_i$ be the contribution to the sum in \eqref{last3} from copies $K$
of $K_k$ such that $K\cap G'$ has exactly $i$ components (including
trivial components of size $1$).
Since $K$ must contain an edge of $I_t\subset E(G')$,
we have $1\le i\le k-1$. Let $\Delta=\Delta(G')\le (\log n)^3$, and set
\[
 z_i = \ell n^{i-1} \Delta^{k-1-i} p^{\binom{k}{2}-\binom{k+1-i}{2}}.
\]
It is easy to check that for $1\le i\le k-1$ we have $Z_i\le z_i$:
there are $\ell$ choices for (one of the) edges of $I_t$ to include,
then picking vertices one by one, either $n$ choices if we start a new
component of $K\cap G'$, or at most $\Delta$ if we do not.
Finally, if $K\cap G'$ has $i$ components, then considering
the case where these components are all complete, by convexity
$E(K\cap G')$ is maximized if the components have sizes $k+1-i,1,1,\ldots,1$,
so $|E(K)\setminus E(G')|\ge \binom{k}{2}-\binom{k+1-i}{2}$.
For $i=1$ we may improve our estimate slightly:
since $G'$ does not contain $K_k$, we gain at least a factor of $p$,
so $Z_1\le pz_1$.

Now $Z_1\le pz_1=p\ell\Delta^{k-2}\le p(\log n)^{O(1)}=o(1)$, so this
contribution to \eqref{last3} is negligible. To handle the remaining cases,
note that
\[
 z_{i+1}/z_i = n\Delta^{-1}p^{k-i},
\]
so $z_{i+1}/z_i$ increases as $i$ increases. Hence the maximum
of $z_i$ for $2\le i\le k-1$ is attained either at $i=2$ or at $i=k-1$.
Now $z_2 =\ell n\Delta^{k-3}p^{k-1}$, and it is easy to check
that this is $o(1)$. Indeed, since $\mu'=\Theta(1)$ we have
$p=\Theta(n^{-\frac{2k}{k(k-1)+2\ell}})$,
and since $2k(k-1)> k(k-1)+2\ell$ we have that $np^{k-1}$ is a constant
negative power of $n$.

At the other end of the range,
$z_{k-1} = \ell n^{k-2}p^{\binom{k}{2}-1}= \Theta(\nu/(n^2p)) = o(1/(np))$,
and we have $np\to\infty$. Thus both $z_2$ and $z_{k-1}$ are $o(1)$,
which gives $z_i=o(1)$ for $2\le i\le k-1$.
Thus
\[
  \sum_K p^{E(K)\setminus (E^+\cup E_t^+\cup I_t)}
 = \sum_{i=1}^{k-1}Z_i \le pz_1+\sum_{i=2}^{k-1}z_i =o(1),
\]
proving \eqref{last3}. This was all that remained to prove the theorem.
\end{proof}

It is natural to wonder whether Theorem~\ref{th_e} can be extended.
For $\ell=1$, the picture is complete: defining $p_0$
by $\mu'(p_0)=1$, since $\mu'=\Theta(\nv p^{\ell})=\Theta(\nv p)$
we have $\nv=\Theta(1/p_0)=o(n)$ whenever $p=\Theta(p_0)$.
As noted earlier, percolation in $\Gokl_p$ for $p/p_0\to\infty$
follows by monotonicity arguments.

For general $k$ and $\ell$, the conditions of Theorem~\ref{th_e}
can presumably be relaxed at least somewhat. Unfortunately, the proof
we have given relies on $\nu=o(n)$, and hence on $\ell< k/2$.

\subsection{Copies of general graphs}

We conclude this paper by briefly considering the graph $G_H^\ell(p)$
obtained from $G(n,p)$ by taking one vertex for each copy of some fixed
graph $H$ with $|H|=k$, and joining two vertices if these copies share at least
$\ell$ vertices, where $1\le \ell\le k-1$.
Ones first guess might be that the results in Section~\ref{sec_cv}
extend at least to regular graphs $H$ without much difficulty, but
this turns out to be very far from the truth. In fact, it seems
that almost all cases are difficult to analyze.

We start with the most interesting end of the range, where $\ell=k-1$,
as in the original question of Der{\'e}nyi, Palla and Vicsek~\cite{DPV}.
To keep things simple, let $H$ be the cycle $C_k$.
For $k=3$, $C_k$ is complete, so this case is covered in Section~\ref{sec_cv}.
The case of $C_4$ is already interesting: when moving from one copy
of $C_4$ to another, we may change opposite vertices essentially independently
of each other. The appropriate exploration is thus as follows:
suppose we have reached a $C_4$ with vertex set $P_0\cup P_1$,
where each of $P_0$ and $P_1$ is a pair of opposite vertices.
Furthermore, suppose we reached this $C_4$ from another $C_4$ containing $P_0$.
Then we continue by replacing
$P_0$ by some other pair $P_2$ of common neighbours of $P_1$.
Suppose that $p=\Theta(n^{-1/2})$; in particular, set $p=\la n^{-1/2}$.
The number $Z$ of common neighbours of $P_1$ outside $P_0$
has essentially a Poisson distribution with mean $\la^2$.
The number of choices for $P_2\ne P_0$ is $\binom{Z+2}{2}-1$, which has expectation
\[
 \E\left(\frac{(Z+2)(Z+1)}{2}-1\right) 
 = \E\left(\frac{Z(Z-1)}{2}+2Z\right) = \la^4/2+2\la^2,
\]
so we believe the critical point will be when this expectation is $1$,
i.e., at
\begin{equation}\label{c4c}
 p=p_0= n^{1/2}\sqrt{\sqrt{6}-2}.
\end{equation}
Of course, it is not clear that the branching process approximation
we have implicitly described is a good approximation to the component
exploration process. However, it is not hard to convince oneself 
that this is the case, at least at first. The key point is that when
we have not yet reached many vertices, the chance of finding a new
vertex adjacent to three or more reached vertices is very small.
Hence the sets of common neighbours of two pairs $P$ and $P'$ are essentially
independent, even if $P$ and $P'$ share a vertex.
We have not checked the details, but we expect that it is not hard
to show rigorously that $p_0$ is indeed the threshold in this case,
although unforeseen complications are of course conceivable.

Taking things further, one might expect the argument above to work for $C_6$,
say, but in fact it breaks down after one step. Suppose we start from
$aubvcw$ and first replace $a$, $b$ and $c$ by other suitable vertices.
Then we have sets $A$, $B$ and $C$ of candidates for $a$, $b$, and $c$.
The problem is at the next step: the possibilities for $u'$, $v'$, $w'$
associated to different triples $(a',b',c')\in A\times B\times C$
are far from independent: for triples $(a',b',c')$ and $(a',b',c'')$,
the choices for $u'$ are exactly the same.
In fact, not only can we not prove a result for any $C_k$, $k\ge 5$, but we
do not even have a conjecture as to the correct critical probability,
although this is clearly of order $\Theta(n^{-1/2})$.

Although $C_4$ is the simplest non-complete example, cycles turn out
not to be the easiest generalization: it is almost certainly not hard
to adapt the outline argument above to complete bipartite
graphs $K_{r,s}$. If $s=r$, then setting $p=\la n^{-1/r}$,
and letting $Z_r$ denote a Poisson
random variable with
mean $\la^r$, the critical point should be given by the solution to
\[
 \E\left(\binom{Z_r+r}{r}-1\right) =1,
\]
generalizing \eqref{c4c}.
If $r\ne s$, the situation is a little different, as alternate
steps in the exploration have different behaviour. Suppose that $r<s$,
and set $p=\la n^{-(s+1)/(rs+s)}$. Then $np^r\to\infty$ and $np^s\to 0$,
so on average a set of $r$ vertices has many common neighbours, 
and so lies in many copies of $K_{r,s}$, while a typical set
of $s$ vertices has no common neighbours.
Starting from a given $K_{r,s}$, with vertex classes $R$ and $S$
of sizes $r=|R|$ and $s=|S|$, let $T$ denote the set of common neighbours
of $R$. Then $\E(|T\setminus S|)=(n-r)p^r\to\infty$, and $|T|$
will be concentrated near $np^r$.
Replacing $S$ by any of the other $\binom{|T|}{s}-1\sim n^sp^{rs}/s!$
subsets $S'$ of $T$ of size $s$, since sets of size $s$ have few common neighbours,
the most likely way the exploration will continue is that some $S'$
will have one common neighbour $x$ outside $R$. Then for each vertex $y$ of $R$
we reach a new $K_{r,s}$ with $R$ replaced by $R\setminus\{y\}\cup\{x\}$.
For each $S'$ we expect to find around $np^s\to 0$ such vertices $x$,
so overall the average number of new choices for $R'$ is
$(1+o(1))n^sp^{rs}s!^{-1} np^s r$, and we expect the critical point
to be given by
$\la n^{-(s+1)/(rs+s)}$ where $\la$ satisfies $\la^{rs+s}=s!/r$;
we have not checked the details.

Finally, since the case $\ell=k-1$ seems too hard in general, one
could consider the other extreme $\ell=1$. This is much easier,
though also less interesting. If $H$ is strictly balanced,
it is very easy to see that the critical point occurs when $(k-1)\mu=1$,
where $\mu$ is the expected number of copies of $H$ containing
a given vertex $v$. For non-balanced $H$ things are a little more complicated:
having found a `cloud' of copies of $H$ containing a single copy of the
(for simplicity unique) densest subgraph $H'$ of $H$,
one next looks for a second cloud meeting the current cloud, and the critical
point should be when the expected number of clouds meeting a given cloud
is $1$. This type of argument can probably be extended to $\ell=2$, at
least if we impose the natural condition in this case that our
copies of $H$ should share an edge, rather than just two vertices.
Beyond this, the whole question seems very difficult.

\end{document}